\documentclass[preprint, 3p]{elsarticle}
\usepackage{setspace}

\usepackage{amsmath}
\usepackage{amsfonts}
\usepackage{bm}
\usepackage[numbers]{natbib}

\usepackage{pgfplots}
\usepackage[lofdepth,lotdepth]{subfig}

\usepackage{tikz}
\usepackage{xcolor}
\usetikzlibrary{shapes,arrows,decorations.pathreplacing}
\definecolor{TUMGreen}{RGB}{162,173,0}
\definecolor{TUMGrey}{RGB}{179,179,179}
\definecolor{TUMDarkGrey}{RGB}{138,138,138}

\usepackage{booktabs}

\usepackage{algorithm,algorithmic} 

\usepackage{hyperref}

\usepackage{float}

\journal{Journal of Computational Physics}

\begin{document}

\begin{frontmatter}

    \title{Sparse Identification of Truncation Errors}

	\author[1]{Stephan Thaler\corref{cor1}}
	\ead{stephan.thaler@tum.de}
	\author[1]{Ludger Paehler\corref{cor1}\fnref{fn1}}
	\ead{ludger.paehler@tum.de}
	\author[1]{Nikolaus A. Adams}
	\ead{nikolaus.adams@tum.de}

	\address[1]{Institute of Aerodynamics and Fluid Mechanics, Technical University of Munich, 85748 Garching, Germany}

	\cortext[cor1]{Shared first author}

	\fntext[fn1]{Corresponding author}

	\begin{abstract}
		This work presents a data-driven approach to the identification of spatial
		and temporal truncation errors for linear and nonlinear discretization schemes
		of Partial Differential Equations (PDEs). Motivated by the central role of
		truncation errors, for example in the creation of implicit Large Eddy schemes, we
 	        introduce the \textit{Sparse Identification of Truncation Errors} (SITE) framework to automatically
		identify the terms of the modified differential equation from simulation data. We build on recent
		advances in the field of data-driven discovery and control of complex systems
		and combine it with classical work on modified differential equation analysis
		of Warming, Hyett, Lerat and Peyret. We augment a sparse regression-rooted approach
		with appropriate preconditioning routines to aid in the identification
		of the individual modified differential equation terms. The construction of such a
		custom algorithm pipeline allows attenuating of multicollinearity effects 
		as well as automatic tuning of the sparse regression hyperparameters using the 
		Bayesian information criterion (BIC). As proof of concept, we constrain 
		the analysis to finite difference schemes and leave other numerical
		schemes open for future inquiry. Test cases include the
		linear advection equation with a forward-time, backward-space discretization,
		the Burgers' equation with a MacCormack predictor-corrector scheme and the Korteweg-de
		Vries equation with a Zabusky and Kruska discretization scheme. Based on variation
		studies, we derive guidelines for the selection of discretization parameters,
		preconditioning approaches and sparse regression algorithms. The results showcase highly
		accurate predictions underlining the promise of SITE for the analysis and
		optimization of discretization schemes, where analytic derivation of modified
		differential equations is infeasible.
	\end{abstract}
	
	\begin{keyword}
		Sparse Regression, Truncation Error, Modified Differential Equation Analysis,
		Data-driven Scientific Computing, Preconditioning
	\end{keyword}
\end{frontmatter}

\section{Introduction}
\label{chap_introduction}

When constructed, modified differential equations (MDEs) provide
valuable insight into the properties of every discretization scheme including
spatial and temporal truncation errors. However, with increasing nonlinearity
of the discretization scheme or the underlying Partial Differential Equation
(PDE), the analytic approach becomes increasingly intractable. Recent advances
in the data-driven learning of differential equations may now allow us to
overcome the drawbacks of modified differential equation analysis (MDEA) by
reformulating the discovery process as symbolic regression. 

MDEA has its roots in von Neumann's stability analysis, which was developed in
the 1940s with its first discussion in O' Brien et al. in 1950
\cite{o1950study}. Realizing the potential of von Neumann's approach, Hirt
proposed a method to connect the stability of nonlinear difference equations
with the form of the truncation error \cite{hirt1968heuristic}. Building on
Hirt's results, Warming and Hyett then established a direct connection between
von Neumann's stability analysis and the symbolic form of the MDE as presented
by Richtmyer and Morton \cite{richtmyer1967difference}, which is the MDEA as we
understand it today. Based on this insight, they showed that the first few
terms of the MDE dominate the properties of the numerical discretization.
Applications of MDEA developed subsequently include increasing accuracy orders
by elimination of leading order truncation error terms
\cite{KLOPFER.06221981,KLOPFER.1983}, enhancing stability using a nonlinear
numerical viscosity term \cite{Majda.1978} and adaptive mesh refinement
\cite{KLOPFER.06221981}. In light of more widespread application of MDEA,
Griffiths and Sanz-Serna examined the limits of MDEA in 1986. They discovered
stability criteria for MDEA and the fundamental insight that by constructing
an MDE we only use a limited amount of information and the MDE can hence not
fully represent the initial PDE discretization \cite{Griffiths.1986}. This
established clear boundaries on the insight which can be derived from MDEA.

In the past decade the development of implicit Large Eddy Simulation (ILES)
turbulence modeling, e.g. Adams et al. \cite{Adams.2004}, led to a renewed
interest in MDEA. Based on the inherent turbulence modeling capability of the
truncation error as shown by Margolin \cite{Margolin.2002}, ILES approaches
tune the discretization scheme to model the subgrid-scale stress-tensor using
the truncation error. At the centre of this approach is the MDEA, but when
considering complex flow configurations, the construction of the MDE becomes
onerous. An automatic construction of the MDE may overcome this drawback.\\

The recent advent of data-driven approaches to the discovery of symbolic
forms may now provide the toolset to construct MDEs for previously intractable
cases in an automatic fashion. Building on breakthrough results of Hod Lipson
and collaborators 
\cite{bongard2005nonlinear,bongard2007automated,schmidt2009distilling} multiple
approaches for the data-driven discovery of PDEs were developed. The core of
our ansatz can be summarized as symbolic regression applied to MDEs. We hence
constrain the introduction to the three approaches which best satisfy these
applicational requirements. These ans\"{a}tze are the sparse regression
approaches of Kutz and Brunton, the physics-informed machine learning of
Karniadakis and Dong's PDE-Net.\\

In 2016, Brunton et al. proposed their initial sparse regression framework
called \textit{Sparse Identification of Nonlinear Dynamics} (SINDy)
\cite{Brunton.2016}, which is a general method for the data-driven
identification of dynamical systems. Using finite differences and polynomial
interpolation to approximate partial derivatives in time and space, this
framework was subsequently extended to PDEs. Said algorithm is named the  
\textit{PDE functional identification of nonlinear dynamics algorithm}
(PDE-FIND) \cite{Rudy.2017} and is explained more closely in section
\ref{chap_theo_PDE-FIND}. Recently, SINDy has seen more generalizations with
extensions to time-evolving parametric PDEs \cite{Rudy.632018} and model
selection using information criteria \cite{Mangan.2017}.

Starting with \textit{Machine Learning of Linear Differential Equations using
Gaussian Processes} \cite{Raissi.2017} in 2017, Raissi, Karniadakis and
collaborators began to develop a diverse array of techniques for the discovery
of coefficients for symbolic terms. This line of inquiry began with a Gaussian
process (GP) construction in which the unknown coefficients are recast as GP
kernel hyperparameters, which can then be learned through optimization of the
marginal likelihood. Notable extensions to this initial framework are the
encoding of time integration schemes in the GP kernel as \textit{Numerical GPs}
\cite{Raissi.2018} and a reformulation of the framework based on Neural
Networks \cite{Raissi.2019}. Raissi et al.'s key insight when working with
neural networks was the use of the symbolic form as an additional loss function
of the network, hence forcing the neural network to obey the physical structure
of the system \cite{Raissi.2019}.

A third distinct approach are the \textit{PDE-Nets} of Long, Lu and Dong
\cite{Long.2018,Long.20181130}. Here, convolutional neural networks are
constructed with partially constrained filters, which approximate the
differential operators. This approach rests on breakthroughs of Cai et
al. \cite{cai2012image} and Dong et al. \cite{dong2017image} in which a direct
connection between filters and finite difference approximations of differential
operators was established. The candidate filters are combined to form the
unknown PDE and then predict the function value at the next time step using
Forward Euler. Learning the filters by minimizing the loss function, we can
then rediscover the exact form of the PDE.\\

The requirements of the task constrain applicable approaches to the sparse
regression framework of Brunton and Kutz \cite{Rudy.2017,Zheng.2019}.
\textit{Numerical GPs} have not seen extension to inverse problems, Neural
Network approaches lack accuracy and the symbolic identification of the
\textit{PDE-Net} cannot be automated. A detailed explanation of the
requirements for the ansatz can be found in section \ref{chap_theo_PDE-FIND}.\\

Building on the results of Rudy et al. \cite{Rudy.2017}, we present a proof of
concept showing that MDEs can be identified from simulation data with high
accuracy and minimal prior knowledge. Such proof of concept is intended as a
stepping stone towards problems for which the analytic derivation of MDEs is
intractable. Our approach, the \textit{Sparse Identification of Truncation
Errors} (SITE) framework could allow for the discovery of the MDE in these
cases. In line with the vision of a recent report by the National Academy of
Sciences \cite{national2012assessing}, we understand the current paper as a
first development towards fully data-driven MDEA tools for the analysis and
optimization of truncation errors. Implicit LES modeling is a direct
application where the proposed algorithm would allow an optimal utilization of
the truncation errors in the construction of the subgrid-scale model. Note that
the outlined procedure is independent of the discretization method. We chose
finite-difference schemes as a framework in which tractable problems can be
defined and analytic MDEs in series form can be derived. Other discretization
methods like the finite element method or the finite volume method can be
analyzed in a similar manner using the SITE approach. For this proof of concept,
we consider one-dimensional test cases only. Extension of SITE to multiple
dimensions is straightforward following the extension of the PDE-FIND algorithm
to multiple dimensions \cite{Rudy.2017}.\\

The paper is structured as follows. Chapter \ref{chap_theory_background} gives
an overview of the theory of analytic MDEA, the preliminaries of the
data-driven identification framework and outlines the major challenges. Chapter
\ref{chap_comp_setup} summarizes the proposed workflow of the SITE approach
before we validate the numerical solvers used for data generation in chapter
\ref{chap_solvers}. Numerical test cases in chapter \ref{Results} demonstrate
the applicability of the procedure to linear and nonlinear PDEs of interest to
fluid dynamics: The advection equation is discretized using a forward-time
backward-space (FTBS) scheme, a MacCormack predictor-corrector scheme
\cite{MACCORMACK.1969} is used for Burgers' equation and the Korteweg-de Vries
equation (KdV) is discretized with the Zabusky and Kruskal scheme
\cite{Zabusky.1965}. Hereby, we assess the applicability of several sparse
regression algorithms to the problem of data-driven identification of MDEs, the
impact of discretization parameters and the effect of preconditioning. In
chapter \ref{chap_Discussion}, we summarize our key findings and give
guidelines for practical applications before presenting an outlook on future
work in chapter \ref{chap_Conclusion}.


\section{Preliminaries}
\label{chap_theory_background}

We will present the required theory, starting with a discussion of the analytic
derivation of MDEs necessary to assess the quality of predictions in our test
cases. This is followed by a coherent exposition of the sparse identification
framework and one of its major challenges, multicollinearity. An overview of
the investigated sparse regression algorithms at the core of the proposed
procedure concludes this preliminaries section.

\subsection{Modified Differential Equation Analysis}
\label{chap_MDEA}

For demonstration purposes, we begin by considering the linear advection
equation
%
\begin{equation}
    u_t + a u_x = 0,
    \label{Advection_equation}
\end{equation}
with a forward-time, backward-space (FTBS) discretization scheme, which is
first order accurate in space and time.
%
\begin{equation}
    \frac{u_i^{j+1} - u_i^j}{\Delta t} + a \frac{u_i^j - u_{i-1}^j}{\Delta x} = 0
    \label{FTBS_Scheme_theory}
\end{equation}
In order to construct the MDEs, existence of a continuously
differentiable function $v(x,t)$ is presumed, which coincides with the
numerical solution obtained from eq. \eqref{FTBS_Scheme_theory} at the
gridpoints $v(x,t) = v(i\Delta x, j \Delta t) = u_{i}^j$ \cite{Warming.1974}.
Substitution in eq. \eqref{FTBS_Scheme_theory} 
%
\begin{equation}
    \frac{v(x, t + \Delta t) - v(x,t)}{ \Delta t} + a \frac{v(x,t)-v(x - \Delta x, t)}{\Delta x} = 0
\end{equation}
and Taylor expansion of each term around $v(x,t)$ yields
\begin{equation}
    \label{first_mod_FTBS}
    v_t + a v_x + v_{tt} \Delta t /2 + v_{ttt} \Delta t^2 /6  + ... - v_{xx} a \Delta x /2 + v_{xxx} a \Delta x^2 /6 + ... = 0,
\end{equation}
dropping the argument of $v(x,t)$. Villatoro et al. \cite{Villatoro.1999} refer
to this form of the MDE in eq. \eqref{first_mod_FTBS} as the first modified
equation, stressing the non-uniqueness of MDEs, since any linear combination of
derivatives of eq. \eqref{first_mod_FTBS} is again a MDE. Following the
procedure of Warming et al. \cite{Warming.1974}, the third modified equation
\cite{Villatoro.1999} can be derived by substitution of higher order time
derivatives with spatial derivatives, using a symbolic mathematics package
(e.g. \textit{sympy}  \cite{Meurer.2017}).
%
\begin{align}
    \begin{split}
        & v_t + a v_x + v_{xx} \Delta x (-a + a^2 h) / 2  + v_{xxx} \Delta x^2 (a - 3 a^2 h + 2 a^3 h^2) / 6 \\
        & + v_{xxxx} \Delta x^3  (-a + 7 a^2 h - 12 a^3 h^2 + 3 a^4 h^3) / 24  \\
        &+ v_{xxxxx} \Delta x^4  (a - 15 a^2 h + 50 a^3 h^2 - 60 a^4 h^3 + 24 a^5 h^4 ) / 120  \\
        &+ v_{xxxxxx} \Delta x^5  (-a^5 + 31 a^2 h - 180 a^3 h^2 + 390 a^4 h^3 - 360 a^5 h^4 + 120 a^6 h^5) / 720  \\
        &+ \mathcal{O}(\Delta x^6) = 0,
        \label{Modified_Advection_equation}
    \end{split}
\end{align}
where $h = \Delta t / \Delta x$. While the procedure of Warming et al.
\cite{Warming.1974} only applies to linear equations, Lerat et al.
\cite{Lerat.1974} previously already derived third MDEs for nonlinear
equations. When truncating the Taylor series in time and space in the first MDE
\eqref{first_mod_FTBS}, higher order initial conditions (ICs) and boundary
conditions (BCs) are necessary in order to obtain a well-posed problem. There
are no higher order time derivatives in a truncated version of the third MDE
\eqref{Modified_Advection_equation} and the second challenge, higher order BCs,
might be handled by enforcing periodic BCs. Yet, there is still criticism
directed at the forward numerical solution of MDEs, because a smooth function
$v(x,t)$ coinciding with the numerical solution on its gridpoints in general
may not satisfy eq. \eqref{Modified_Advection_equation} \cite{Chang.1990}.
However, one is usually not interested in solutions of the MDE for specific ICs
and BCs, but rather in the form of the MDE itself that contains information
about the discretization scheme \cite{Chang.1990}. Considering this, we
demonstrate that reconstruction of MDEs from data is possible without referring
to any specific forward solution of the respective MDE.

\subsection{Identification Framework}
\label{chap_theo_PDE-FIND}

The properties of the problem at hand largely determine the choice of an
appropriate symbolic identification framework. The framework has to be
applicable to nonlinear equations, deal with a moderate to large number of
candidate functions and yield high accuracy predictions in a low noise
environment. High accuracy is essential, given that coefficients in MDEs often
span several orders of magnitude, e.g. table \ref{tab:Advection_model}. GPs of
Raissi et al. \cite{Raissi.2017} do not allow nonlinear equations,
\textit{Numerical GPs} \cite{Raissi.2018} are not yet extended to inverse
problems and \textit{Physics Informed Neural Networks'} (PINNs)
\cite{Raissi.2019} continuous and discrete time models suffer from the accuracy
requirement. Using PINNs, we only managed to discover terms up to $v_{xxx}$ in
eq. \eqref{Modified_Advection_equation}, all higher order derivatives were not
identified correctly. We attribute this behavior to the large number of
trainable parameters in the neural network, since finding the precise optimum
in high-dimensional spaces is a difficult task. Including high order
derivatives aggravates the issue, because the effective depth of the network
increases exponentially with increasing derivative order, compromising neural
network training. 

For our work, we choose PDE-FIND \cite{Rudy.2017}. With the numerical test
cases in chapter \ref{Results}, we will demonstrate that all of the outlined
requirements are met. The algorithm leverages linear sparse regression, which
facilitates accurate parameter estimation and sets terms not included in the
predicted sparse equation exactly to 0. The design of the library enables
flexible construction of fairly general candidate function spaces, including
higher order temporal derivatives for identification of first MDEs.
Furthermore, the vast amount of literature on (sparse) linear regression
provides a firm theoretical basis for our work. \\ 

\tikzstyle{sensor}=[draw, fill=blue!20, text width=5em, 
    text centered, minimum height=2.5em]
\tikzstyle{ann} = [above, text width=5em]
\tikzstyle{ContentColumn} = [sensor, text width=1em, fill=white, draw=TUMDarkGrey,
    minimum height=18em, rounded corners]
\tikzstyle{DotDotColumn} = [sensor, text width=1em, fill=white, draw=white,
	minimum height=18em, rounded corners]
\def\blockdist{0}
\def\edgedist{0}

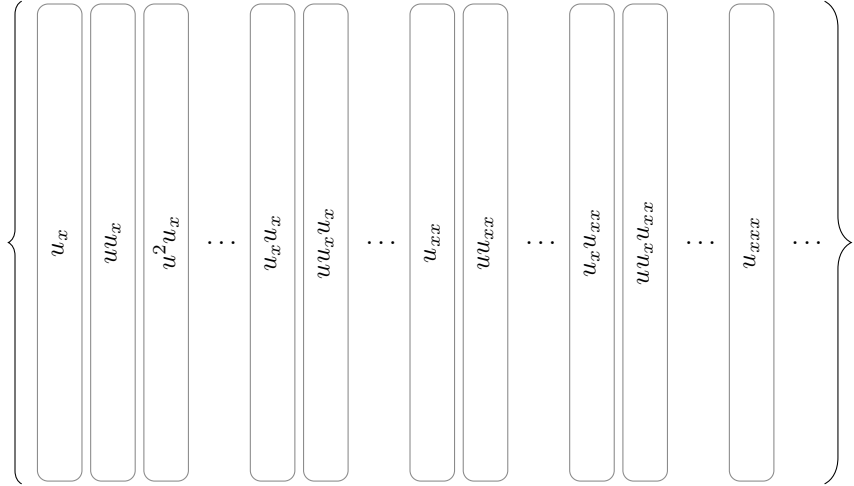
\begin{figure}[H]
	\centering
	\begin{tikzpicture}[node distance=0.7cm]
    	\node [ContentColumn] (content1) {\rotatebox{90}{$u_{x}$}}; 
		\node [ContentColumn, right of=content1] (content2) {\rotatebox{90}{$u u_{x}$}};
		\node [ContentColumn, right of=content2] (content3) {\rotatebox{90}{$u^{2}u_{x}$}};
		\node[DotDotColumn, right of=content3] (DotDot1) {\dots};
    	\node [ContentColumn, right of=DotDot1] (content4) {\rotatebox{90}{$u_{x} u_{x}$}};
    	\node [ContentColumn, right of=content4] (content5) {\rotatebox{90}{$u u_{x} u_{x}$}};
    	\node[DotDotColumn, right of=content5] (DotDot2) {\dots};
		\node[ContentColumn, right of=DotDot2] (content6) {\rotatebox{90}{$u_{xx}$}};
		\node[ContentColumn, right of=content6] (content7) {\rotatebox{90}{$u u_{xx}$}};
		\node[DotDotColumn, right of=content7] (DotDot3) {\dots};
		\node[ContentColumn, right of=DotDot3] (content8) {\rotatebox{90}{$u_{x}u_{xx}$}};
		\node[ContentColumn, right of=content8] (content9) {\rotatebox{90}{$u u_{x} u_{xx}$}};
		\node[DotDotColumn, right of=content9] (DotDot4) {\dots};
        \node[ContentColumn, right of=DotDot4] (content10) {\rotatebox{90}{$u_{xxx}$}};
		\node[DotDotColumn, right of=content10] (DotDot4) {\dots};
    	\draw [decorate,decoration={brace,amplitude=5pt}]
    	(-0.5,-3.2) -- (-0.5,3.2) node [black,midway,xshift=-0.6cm] 
    	{};
    	\draw [decorate,decoration={brace,amplitude=10pt,mirror},xshift=-4pt,yshift=0pt]
    	(10.2,-3.2) -- (10.2,3.2) node [black,midway,xshift=-0.6cm] 
    	{};
	\end{tikzpicture}
    \caption{Example library $\bm\Theta(\mathbf{u}) \in \mathbb{R}^{n \times p}$.
             In our test cases, the library contains from 49 to 210 candidate
             terms for the advection equation, 28 terms for the Burgers'
             equation and 68 terms for the KdV equation.}
    \label{fig:Example_Library}
\end{figure}

We summarize the PDE-FIND method, following \cite{Rudy.2017}. Our objective is
to identify first or third MDEs\footnote{Beware non-uniqueness issues if
multiple forms of MDEs (e.g. first and third) can be represented by
$\bm\Theta(\mathbf{u})$.} of the form 
%
\begin{equation}
    u_t = F(u_{tt}, u_{ttt}, ..., u, u_x, u^k u_{xx}, ...). 
    \label{PDE_general_equation}
\end{equation}
Given data $u(x,t)$ obtained from the numerical solver under investigation, eq.
\eqref{PDE_general_equation} is set up at each point on the space-time grid
using finite differences to approximate any partial derivatives in space and
time. This builds a linear system of equations,
%
\begin{equation}
    \label{SINDy System of equations}
    \mathbf{u_t} = \bm\Theta(\mathbf{u}) \bm\xi + \bm\epsilon \  ; \qquad  \bm \epsilon \sim \mathcal{N}(\mathbf{0},\,\sigma^{2} \mathbf{I}),
\end{equation}
where $\mathbf{u} \in \mathbb{R}^n$ is a discretized version of $u(x,t)$. This
system is solved for the unknown weights of candidate terms
$\bm\xi \in \mathbb{R}^n$ using sparse regression methods, given that the
majority of candidate therms in the library
$\bm\Theta(\mathbf{u}) \in \mathbb{R}^{n \times p}$ are not part of the correct
MDE. Fig. \ref{fig:Example_Library} illustrates an example
library $\bm\Theta(\mathbf{u})$. The noise $\bm\epsilon$ is assumed to be
Gaussian; note that this is an implicit assumption in most sparse regression
algorithms. $\bm\epsilon$ can be attributed to higher order terms of the MDE
not  included in the library as well as truncation and round-off error of the
finite difference approximations. Correlated candidate feature  vectors have
been identified as a major challenge for PDE-FIND \cite{Rudy.2017}.
Multicollinearity is a central concern in the construction of
$\bm\Theta(\mathbf{u})$, since its amount increases with the size of the
candidate term pool.

\subsection{Multicollinearity}
\label{sec:multicolinearity}
Multicollinearity has been extensively studied in the
context of linear regression. For example, it is well known that ordinary least
squares (OLS) regression accuracy suffers from it. In-depth coverage of the
topic is widely available; hence we will only provide some exemplary intuition
about the ill-posed nature of multicollinearity:

The sparse regression task in eq. \eqref{SINDy System of equations} is to find
the best approximation to $u_t$ by weighted spatial derivatives of $u(x,t)$.
Consider a particularly weak IC choice $\overline{u}(x,0)=\sin(kx)$. Given two
collinear terms $\overline{u}_{xx} = -k^2 \overline{u}$, any $\gamma \in [0,1]$
in the linear combination of both terms
$\gamma * \overline{u} - (1-\gamma)\overline{u}_{xx}/k^2$ yields the exact same
result. It is thus impossible to distinguish the importance of one term versus
the other. This intuition extends in a straightforward manner to
multicollinearity, where a given feature vector can be represented with little
error by a linear combination of the other feature vectors. The variance
inflation factor (VIF) quantifies multicollinearity by measuring the deviation
of this representation. A large VIF indicates a small deviation and therefore
high multicollinearity.
%
\begin{equation}
    \label{VIF_eq}
    VIF_i = 1-1/(1-R_i^2), 
\end{equation}
where $R_i^2$ is the coefficient of determination of the linear regression
%
\begin{equation}
    \mathbf{f}_i = \xi_{i0} + \sum_{j=1; j \neq i}^p \xi_{ij} \mathbf{f}_j,
\end{equation}
and $\mathbf{f}_i$ is a column of $\bm \Theta(\mathbf{u})$.\\

The most prominent approach
to obtain a well-posed problem is $L_2$ regularization, i.e. ridge regression
\cite{Hoerl.1970}. Generalizing $L_2$ regularization, $L_q$ regularization of
the least-squares minimization problem of eq. \eqref{SINDy System of equations}
is defined by adding a penalty term for the weight vector.
%
\begin{equation}
    \label{Lq penalty unconstrained}
    \underset{\bm\xi}{\arg \min} \left( \frac{1}{2} ||\bm \Theta(\mathbf{u}) \bm\xi - \mathbf{u}_t||_2^2 + \lambda \sum_{i=1}^p |\xi_i|^q \right) \ ; \quad q\geq 0
\end{equation}
$L_2$ penalizes large coefficients (note the $-k^2$ scaling in the example),
resolving the ill-posedness. The solution of eq.
\eqref{Lq penalty unconstrained} can be interpreted from a Bayesian statistics
viewpoint to be  the maximum a posteriori estimate of a Gaussian likelihood and
a log-prior distribution $\log(p(\bm \xi)) = \lambda ||\bm\xi||^q$
\cite[p.~72]{Hastie.2009}. The $L_2$ penalty corresponds to a Gaussian prior on
weights centered at 0, with scale defined by $\lambda$. The prior mean of 0 is
not supported by actual prior knowledge, thus artificially biasing non-zero
coefficients towards zero.

A proper choice of prior for the MDE
identification problem is challenging, since information about the scale, and
ideally the mean, of each coefficient is necessary. The order of the
discretization scheme provides some information about the scale of the dominant
order truncation error terms. However, this value is in general too large for
sufficient regularization of higher order terms, given that coefficients in
MDEs often span several orders of magnitude.

\subsection{Sparse Regression Algorithms}
\label{Theo Sparse Regression} The aim of sparse regression is to find solutions
to eq. \eqref{SINDy System of equations}, such that the fewest features
possible are included in the model, while representing the data
as faultless as achievable given candidate library $\bm \Theta(\mathbf{u})$. An
important criterion here is \textit{sign consistency}, which quantifies the
ability of the algorithm to distinguish relevant from irrelevant features in the
limit $n \rightarrow \infty$ \cite{Zhao.2006}. To promote sparsity, $L_0$
regularization \eqref{Lq penalty unconstrained} represents a natural choice as
it penalizes non-zero elements in $\bm \xi$. $L_0$ regularization can be
formulated either as an unconstrained \eqref{L0 penalty unconstrained} or
constrained \eqref{L0 penalty constrained} optimization problem:
%
\begin{equation}
    \label{L0 penalty unconstrained}
    \underset{\bm\xi}{\arg \min} \left( \frac{1}{2} ||\bm \Theta(\mathbf{u}) \bm\xi - \mathbf{u}_t||_2^2 + \lambda ||\bm\xi||_0 \right) ,
\end{equation}
%
\begin{equation}
    \label{L0 penalty constrained}
    \underset{\bm\xi}{\arg \min} \left( \frac{1}{2} ||\bm \Theta(\mathbf{u}) \bm\xi - \mathbf{u}_t||_2^2 \right) \mathrm{, \ subject \  to} \   ||\bm \xi||_0 \leq s_0.
\end{equation}
Using a leaps and bound algorithm \cite{Furnival.1974}, solution of eq.
\eqref{L0 penalty constrained} is possible up to $s=30\ ..\ 40$
\cite[p.~57]{Hastie.2009}. However, due to being a combinatorically large
problem, most modern sparse regression algorithms focus on numerically
efficient approximations to the $L_0$ regularized problem \cite{Zhang.2009}.
We compare the performance of four sparse regression algorithms in the
identification of MDEs task. 

The most widespread algorithm among the considered is the \textit{Least
Absolute Shrinkage and Selection Operator} (Lasso) \cite{Tibshirani.1996},
which uses $L_1$ regularization. This is the smallest $q$ such that eq.
\eqref{Lq penalty unconstrained} is a convex problem, improving numerical
efficiency of the optimization over non-convex problems significantly. Unlike
ridge regression, Lasso sets irrelevant coefficients exactly to 0
\cite[p.~73]{Hastie.2009}. Lasso has two major deficiencies though. The $L_1$
penalty can be interpreted as a Laplace distribution centered at 0 and its
scale defined by $\lambda$. However, analogous to ridge regression, this prior
is not supported by actual prior knowledge. If the prior mean deviates from the
data, which it does for all non-zero coefficients, this introduces bias that
can be unsatisfactorily large for the problems considered in this paper (see
e.g. fig. \ref{fig:accuracy_advection}). Besides, Lasso is sensitive to
correlated feature vectors: A necessary and almost sufficient condition for
\textit{sign consistency} is the \textit{irrepresentability condition}, which
$\bm \Theta(\mathbf{u})$ has to fulfill. The \textit{irrepresentability
condition} has several sufficient conditions, a common property of which is a
bound on the maximum correlation between feature vectors of
$\bm \Theta(\mathbf{u})$ \cite{Zhao.2006}.

To cope with highly correlated feature vectors common to
$\bm \Theta(\mathbf{u})$, Brunton et al. \cite{Brunton.2016,Rudy.2017} propose
\textit{Sequential Threshold Ridge Regression} (STRidge). Hereby, the $L_2$
penalty is used for regularization of correlated feature vectors, while
sequential thresholding promotes sparsity. Given a tolerance \textit{tol} and
$\lambda$ corresponding to the $L_2$ penalty as hyperparameters, analytically
tractable ridge estimates are obtained sequentially and  weights smaller than
\textit{tol} are set to 0, until the non-zero weights converge and yield a
sparsity pattern. In a second step, an OLS estimate for the non-zero weights is
calculated, thus avoiding any bias from the $L_2$ penalty. 

The recently published \textit{Sparse Relaxed Regularized Regression} (SR3)
method \cite{Zheng.2019} relaxes regularization penalties to decouple the
accuracy and sparsity requirements by introducing an auxiliary weight vector
$\mathbf{w}$. Application of SR3 to problem \eqref{Lq penalty unconstrained}
yields
%
\begin{equation}
    \label{SR3_raw}
    \underset{\bm\xi, \mathbf{w}}{\arg \min} \left( \frac{1}{2} ||\bm \Theta(\mathbf{u}) \bm\xi - \mathbf{u}_t||_2^2 + \lambda \sum_{i=1}^p |w_i|^q  + (\gamma/2) ||\bm \xi -\mathbf{w}||_2^2 \right) \ ; \quad q\geq 0
\end{equation}
where the additional hyperparameter $\gamma$ controls the amount of deviation
of $\bm \xi$, which enforces accuracy, from $\mathbf{w}$, which enforces
sparsity. Minimization of eq. \eqref{SR3_raw} with respect to $\bm \xi$ yields
%
\begin{equation}
    \label{SR3_minimized}
    \underset{\mathbf{w}}{\arg \min} \left( \frac{1}{2} ||\mathbf{F}_{\gamma} w - \mathbf{g}_{\gamma}||_2^2 + \lambda \sum_{i=1}^p |w_i|^q  + (\gamma/2)\right) \ ; \quad q\geq 0
\end{equation}
with 
\begin{equation}
    \kappa(\mathbf{F}_{\gamma}) = \kappa(\bm \Theta(\mathbf{u})) \sqrt{\frac{\gamma + \sigma_{\mathrm{min}}(\bm \Theta(\mathbf{u}))}{\gamma + \sigma_{\mathrm{max}}(\bm \Theta(\mathbf{u}))}},
\end{equation}
where $\kappa = \sigma_{\mathrm{max}} / \sigma_{\mathrm{min}}$ is the condition
number. $\gamma$ controls the amount of reduction in $\kappa$. For the
definitions of $\mathbf{F}_{\gamma} \in \mathbb{R}^{n \times p}$ and
$\mathbf{g}_{\gamma} \in \mathbb{R}^n$ see Zheng et al. \cite{Zheng.2019}. SR3 has
remarkable similarities to the puffer transformation in section
\ref{chap_puffer}, transforming both $\bm \Theta(\mathbf{u})$ and
$\mathbf{u}_t$ and significantly reducing $\kappa$ in the transformed system. 

The Forward-Backward Greedy Algorithm (FoBa) \cite{Zhang.2009} approximates
problem \eqref{L0 penalty constrained} by sequential greedy selection of the
feature, which reduces the OLS residual most. The rationale is to obtain a
sparse solution by maximizing the gain from each added new candidate. After a
fixed number of forward steps, backward steps aim to eliminate the least
important features obtained from forward steps. Features are only deleted if
the residual increases by less than half of the residual decrease of the last
forward steps, guaranteeing convergence in a finite number of steps. This
overcomes one major weakness of forward greedy algorithms, namely the inability
to delete features deemed irrelevant during forward stepping. See Zhang
\cite{Zhang.2009} for an example of selected features deemed irrelevant. FoBa
terminates once the next best forward step reduces the residual by less than
the hyperparameter $\epsilon$. For \textit{sign consistency}, FoBa needs to
meet the \textit{sparse eigenvalue condition} \cite{Zhang.2009}, which is a
weaker assumption than the \textit{irrepresentability condition}
\cite{Bickel.2009}.


\section{Sparse Identification of Truncation Errors}
\label{chap_comp_setup}

The centerpiece of our SITE 
approach is preconditioning at multiple stages of the workflow illustrated in
fig. \ref{fig:AlgoSetup}. For detailed explanations of the individual steps we
refer to the respective subsections in this chapter. We begin by using a
\textit{Non-Uniform Rational Basis Spline} (NURBS) as IC, which has been
constructed to minimize multicollinearity. This is followed by the solver
forward run generating the data. Afterwards, we assemble our candidate library
following the PDE-FIND framework, scale the candidate terms and apply a puffer
transformation. The preconditioned system is solved with a sparse regression
algorithm for various hyperparameter values, resulting in a set of models that
encode the predicted MDEs. From this set we choose the best model based on the
well known Bayesian information criterion (BIC) \cite{Schwarz.1978}.

\tikzstyle{sensor}=[draw, fill=blue!20, text width=5em, 
    text centered, minimum height=2.5em]
\tikzstyle{ann} = [above, text width=5em]
\tikzstyle{ContentColumn} = [sensor, text width=1.2em, fill=white, draw=TUMDarkGrey,
    minimum height=18em, rounded corners]
\tikzstyle{DotDotColumn} = [sensor, text width=1.2em, fill=white, draw=white,
    minimum height=18em, rounded corners]
\tikzstyle{DotDotColumn2} = [sensor, text width=9em, fill=white, draw=white,
    minimum height=9em]
\def\blockdist{2}
\def\edgedist{2}

\begin{figure}[t]
    \centering
	\begin{tikzpicture}[node distance=1cm]
        \node[DotDotColumn2] (Solver) {\scalebox{0.628}{\input{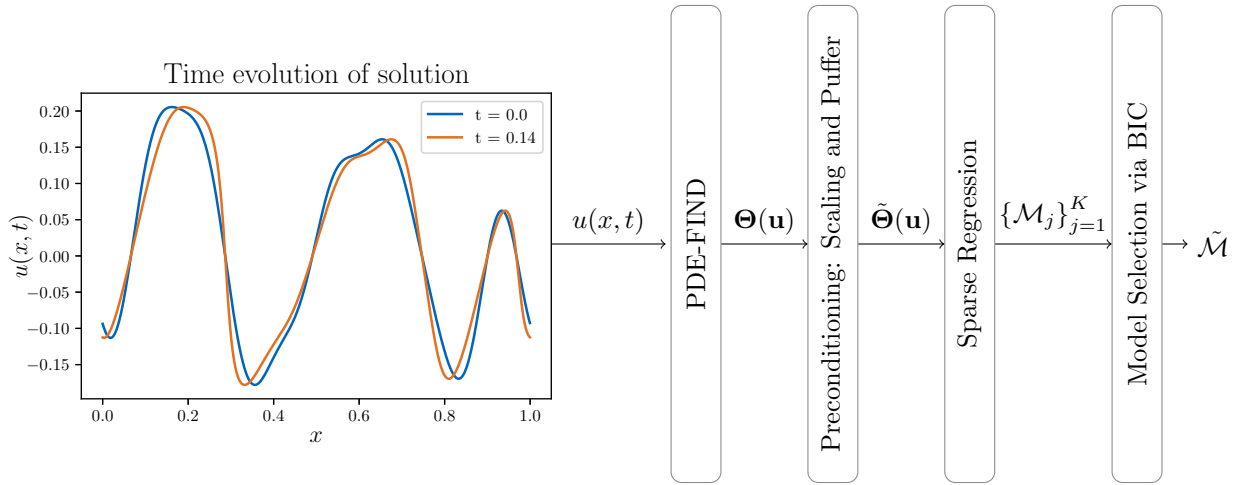}}};
        \node[ContentColumn, right of=Solver, node distance=7.5cm] (PDEFind) {\rotatebox{90}{PDE-FIND}};
        \node[ContentColumn, right of=PDEFind, node distance=1.8cm] (Preconditioning) {\rotatebox{90}{Preconditioning: Scaling and Puffer}};
        \node[ContentColumn, right of=Preconditioning, node distance=1.8cm] (Sparse) {\rotatebox{90}{Sparse Regression}};
        \node[ContentColumn, right of=Sparse, node distance=2.2cm] (Model) {\rotatebox{90}{Model Selection via BIC}};
        \node[DotDotColumn, right of=Model, node distance=1cm] (BestModel) {$\tilde{\mathcal{M}}$};
        \draw[->] (5.6,0) -- (7.1,0) node[midway,above] {$u(x,t)$};
    	\path[->] (PDEFind) edge node[midway,above] {$\bm \Theta (\mathbf{u})$} (Preconditioning);
    	\path[->] (Preconditioning) edge node[midway,above] {$ \tilde{\bm \Theta}(\mathbf{u})$} (Sparse);
    	\path[->] (Sparse) edge node[midway,above] {$\left\{ \mathcal{M}_j \right\}^{K}_{j=1}$} (Model);
    	\path[->] (Model) edge node {} (BestModel);
	\end{tikzpicture}
    \caption{Illustration of the SITE approach, starting with the solver output
             data and resulting in the selected model encoding the predicted MDE}
	\label{fig:AlgoSetup}
\end{figure}

\subsection{Preconditioning}

Given that we cannot rely on prior knowledge to deal with this ill-conditioned
regression problem \eqref{SINDy System of equations}, we propose a three step 
preconditioning procedure to reduce multicollinearity.

\subsubsection{Spline Initialization}
A function $u(x,t)$ that cannot be well described by linear combinations of its
own derivatives is advantageous for low multicollinearity (section
\ref{sec:multicolinearity}). Since data from few simulation time steps are
sufficient for SITE (see chapter \ref{Results}): $u(x,0) \approx u(x,T)$. In
cases in which the MDE does not depend on the current function value
$u(x,t_{i-1})$, the initial condition $u(x,0)$ can be chosen freely. This is an
opportunity to choose $u(x,0)$, such that a well-conditioned problem is
obtained. We use root mean square-VIF (RMS-VIF) \eqref{eq:RMS_VIF} as an
objective function to optimize the IC.
%
\begin{equation}
    \mathrm{RMS-VIF} = \sqrt{\frac{1}{p} \sum_{i=1}^p \mathrm{VIF}_i^2}
    \label{eq:RMS_VIF}
\end{equation}
We employ a NURBS $s_\phi$ as a parametrization for $u(x,0)$, which provides a
flexible definition for $u(x,0)$ by variation of its weight vector $\phi$. This
allows an unconstrained optimization, because high order differentiability is
achieved through the order of the spline and periodicity can be enforced
automatically via additional knots outside $\Omega$, due to local support of
the NURBS. The choice of $u(x,0)$ therefore reduces to finding a weight vector
$\bm \phi$ that minimizes RMS-VIF. We apply a gradient-free particle swarm
optimization algorithm \cite{Zeugmann.2011}, as obtaining the gradient of
RMS-VIF with respect to $\bm \phi$ is non-trivial. The procedure is outlined in
algorithm 1. Given $\bm \phi$, the parametric NURBS $s_\phi$ is created from
which $u(x,0)$ is obtained. Next, the simulation is run and
$\bm\Theta(\mathbf{u})$ is assembled. Lastly, RMS-VIF of
$\bm\Theta(\mathbf{u})$ is calculated and fed back into the particle swarm
optimization, yielding new proposals for $\bm \phi$. Note that the reduction in
RMS-VIF and condition number achieved through NURBS initialization comes at no
additional cost other than increased preprocessing time.

\begin{algorithm}
    \caption{NURBS initialization}
    \begin{algorithmic}[1]
        \renewcommand{\algorithmicrequire}{\textbf{Input:}}
        \renewcommand{\algorithmicensure}{\textbf{Output:}}
        \REQUIRE iters, particles
        \ENSURE  $u(x,0)$
        \\ \textit{Initialization} : $\{\bm \phi^1_j\}_{j=1}^{\mathrm{particles}} = \mathrm{random}$
        \FOR {$i = 1$ .. iters}
        \FOR {$j = 1$ .. particles}
        \STATE $s_\phi$ = NURBS$(\bm \phi_j^i)$
        \STATE $u(x,0)$ = interpolate$(s_\phi)$
        \STATE $u(x,t)$ = forwardsolve($u(x,0)$)
        \STATE $\bm \Theta(\mathbf{u})$ = PDE-FIND$(u(x,t))$
        \STATE RMS-$\mathrm{VIF}^i_j$ = calculateVIF$(\bm \Theta(\mathbf{u}))$
        \ENDFOR
        \STATE $\{\bm \phi^{i+1}_j\}_{j=1}^{\mathrm{particles}}$ = particleswarm$(\{$RMS-$\mathrm{VIF}^i_j\}_{j=1}^{\mathrm{particles}})$
        \ENDFOR
        \STATE $\bm \phi^{\mathrm{best}}$ = $ \underset{j}{\arg \min} $ RMS-$\mathrm{VIF}(\{\bm \phi_j^{\mathrm{iters}}\}_{j=1}^{\mathrm{particles}})$
        \STATE $s_\phi$ = NURBS$(\bm \phi^{\mathrm{best}})$
        \STATE $u(x,0)$ = interpolate$(s_\phi)$
        \RETURN $u(x,0)$
    \end{algorithmic} 
\end{algorithm}

\subsubsection{Scaling}

Higher order derivatives of $u(x,t)$ often span several orders of magnitude.
Scaling $\bm\Theta(\mathbf{u})$ such that all features share a common
magnitude is a standard preprocessing step in regression analysis
\cite{Neumaier.1998,Hastie.2009}. We can rewrite eq.
\eqref{SINDy System of equations} using a scaling matrix
$\mathbf{S} \in \mathbb{R}^{p \times p}$ 
%
\begin{equation}
    \label{scaled equation}
    \mathbf{u}_t = \overline{\bm \Theta}(\mathbf{u}) \overline{\bm\xi} + \bm \epsilon \  ; \quad  \bm \epsilon \sim \mathcal{N}(\mathbf{0},\,\sigma^{2} \mathbf{I}) \quad \mathrm{with} \quad \overline{\bm \Theta}(\mathbf{u}) = \bm \Theta(\mathbf{u}) \mathbf{S}^{-1} \  ; \quad \overline{\bm\xi} = \mathbf{S} \bm \xi,
\end{equation}
where we use a default diagonal scaling matrix \cite{Neumaier.1998}
%
\begin{equation}
    \mathbf{S}_{kk} = \sqrt{(\bm \Theta(\mathbf{u})^T \bm \Theta(\mathbf{u}))_{kk}}.
\end{equation}
Even though rescaling does not reduce VIF, it reduces the condition number
$\kappa(\overline{\bm\Theta}(\mathbf{u})) \leq \kappa(\bm\Theta(\mathbf{u}))$.
If $\kappa(\bm\Theta(\mathbf{u}))$ is very large, finding a numerical solution
of OLS is problematic \cite{Golub.1999,Paige.1982}.

\subsubsection{Puffer Transformation}

\label{chap_puffer} Puffer transformation \cite{Jia.8282012} is one of our main 
preconditioning steps. While originally being introduced to extend the applicability
of Lasso to problems, where the \textit{irrepresentability condition} is
not met, we show that other sparse regression algorithms benefit from the
puffer transformation as well (see fig. \ref{fig:accuracy_advection} ). Eq.
\eqref{scaled equation} is rewritten by multiplying a precondition matrix
$\mathbf{F} \in \mathbb{R}^{n \times n}$ from the left
%
\begin{equation}
    \label{norm_precondition}
    \tilde{\mathbf{u}}_t = \tilde{\bm \Theta}(\mathbf{u}) \overline{\bm \xi} + \tilde{\bm \epsilon} \ ; \quad \mathrm{with} \quad \tilde{\mathbf{u}}_t = \mathbf{F} \mathbf{u}_t \  ; \  \tilde{\bm \Theta}(\mathbf{u}) = \mathbf{F} \overline{\bm \Theta}(\mathbf{u}) \  ; \  \tilde{\bm \epsilon} = \mathbf{F} \bm \epsilon \sim \mathcal{N}(\mathbf{0},\,\tilde{\bm \Sigma}).
\end{equation}
Jia et al. \cite{Jia.8282012} propose to construct $\mathbf{F}$ from a singular
value decomposition (SVD)
$\overline{\bm \Theta}(\mathbf{u}) = \mathbf{U} \mathbf{D} \mathbf{V}^T$;
$\mathbf{U} \in \mathbb{R}^{n \times p}; \mathbf{D} \in \mathbb{R}^{p \times p}; \mathbf{V}^T \in \mathbb{R}^{p \times p}$,
assuming $\mathrm{rank}(\overline{\bm \Theta}(\mathbf{u}))  = p$ with $n > p$
and neglecting zero rows and columns \footnote{If $n < p$, the design matrix is
projected onto the Stiefel manifold, resulting in an empirically verifiable
significant reduction in pairwise feature correlation \cite{Jia.8282012}.}
%
\begin{equation}
    \label{F_preconditioner}
    \mathbf{F} = \mathbf{U} \mathbf{D}^{-1} \mathbf{U}^T  \  \Longrightarrow \tilde{\bm \Sigma} = \sigma^2 \mathbf{U} \mathbf{D}^{-2} \mathbf{U}^T.
\end{equation}
Given that $\mathbf{U}^T \mathbf{U} = \mathbf{V}^T \mathbf{V} = \mathbf{I}_d$
from the SVD, it is straightforward to show orthonormality of the
preconditioned system matrix
$\tilde{\bm \Theta}^T(\mathbf{u}) \tilde{\bm \Theta}(\mathbf{u}) = \mathbf{I}_d$,
yielding a perfect $\kappa(\tilde{\bm \Theta}(\mathbf{u}))=1$ (neglecting
numerical error). However, the elimination of multicollinearity  comes at the
cost of inflating the noise variance $\tilde{\bm \Sigma}$ by $\mathbf{D}^{-2}$.
This reduces the signal-to-noise ratio, thus counteracting the benefits from a
better conditioned problem \cite{Jia.8282012}. At this point, the two previous 
preconditioning steps become relevant. Reducing multicollinearity by a good
choice of $u(x,0)$ and scaling $\bm \Theta(\mathbf{u}))$ both reduce
$\kappa(\overline{\bm \Theta}(\mathbf{u}))$. Therefore, the smallest singular
values in $\mathbf{D}$ become larger, hence reducing the noise inflation effect.

The preconditioner $\mathbf{F}$ is generalized by Jia et al. \cite{Jia.2015},
with the goal to bound the inflation of the noise caused by too small singular
values: $\mathbf{\hat{D}}^{-1}$ substitutes $\mathbf{D}^{-1}$ in eq.
\eqref{F_preconditioner},
%
\begin{equation}
    \mathbf{\hat{D}}^{-1}_{kk} = g(\mathbf{D}_{kk},\tau) / \mathbf{D}_{kk},
\end{equation}
requiring a reasonable choice for $g$ and $\tau$. While this approach might be
a helpful remedy if noise inflation is a major issue, for the sake of
simplicity, we restrict ourselves in this work to the default $\mathbf{F}$ in
eq. \eqref{F_preconditioner}.

\subsection{Sparse Regression}

The purpose of the sparse regression algorithm is to propose a sparse solution
to the linear system  \eqref{SINDy System of equations} $\bm \xi$ encoding the
predicted MDE. We chose FoBa as the default sparse regression algorithm for
SITE due to its high accuracy, numerical efficiency and straightforward
hyperparameter tuning (to be shown in chapter \ref{Results}). Iterating over
the hyperparameter $\epsilon$ of FoBa gives a set of candidate models
$\{\mathcal{M}_j\}_{j=1}^K$, with a different degree of sparsity  each. Since
the number of terms to be included in the true model\footnote{We agree that a
model can never perfectly represent reality and therefore a "true" model
usually does not exist. Since for the problems considered in this paper we can
derive the exact solution analytically, we are in fact dealing with the rare
case that a true model exists (yet having an infinite number of parameters). We
will refer to "true" terms for those terms included in the analytic model.} is
not known a priori, we employ a data-driven model selection procedure.

\subsection{Model Selection}
We calculate the BIC for all linear models $\{\mathcal{M}_j\}_{j=1}^K$ with
respect to a single test simulation and select the model that maximizes the
BIC. A test dataset $\tilde{\bm \Theta}^{\mathrm{test}}(\mathbf{u})$ can be
computed easily  by changing $u(x,0)$. In this paper, we use a NURBS with a
smaller number of optimizable knots $\phi$ resulting in a different initial
condition. For a linear regression model, BIC is defined as
\cite[p.~153]{Berger.2001}
%
\begin{equation}
\label{eq_BIC}
    \mathrm{BIC} = -\frac{n_{\mathrm{eff}}}{2} \log \left( |\tilde{\bm \Theta}^{\mathrm{test}}(\mathbf{u}) \hat{\bm \xi} - \tilde{\mathbf{u}}_t^{\mathrm{test}} |_2^2 \right) - \frac{k}{2} \log \left( n_{\mathrm{eff}} \right),
\end{equation}
where $n_{\mathrm{eff}}$ is the effective sample size, $k$ is the number of
features included in the model and $\hat{\bm \xi}$ is the sparse regression
result from the training set. The BIC aims to find parsimonious models from a
balance of the residual sum of squares in the first term and model complexity
in the second term of eq. \eqref{eq_BIC}. Only for independent, identically
distributed data $n_{\mathrm{eff}} = n$ \cite{Berger.2014}. Clearly, the
number of independent samples is smaller than $n$: Considering the case of
ideal advection and periodic boundaries, $u(x,0)$ is only shifted over time.
Therefore, all samples of future time steps are identical to the samples from
the first time step (with the simplifying assumption that the advection
distance is a multiple of $\Delta x$). Due to nonlinear equations and higher
order terms in the MDE, $u(x,t)$ slightly deforms over time, thus increasing
the amount of information compared to the ideal advection case. However,
correlation between samples within one time step (e.g. neighboring points)
decreases the amount of information contained in the data from a single time
step. Presuming these two secondary effects roughly balance each other, we
take $n_{\mathrm{eff}} = n^{x}$ as a coarse estimate of $n_{\mathrm{eff}}$,
where $n^{x}$ is the spatial resolution of $\Omega$.\\

If the puffer transformation is applied, the noise vector
$\mathbf{F}\bm\epsilon$ contains statistically dependent terms, which are not
accounted for in the assumptions of BIC. Jia et al. \cite{Jia.8282012}
therefore propose to compute an OLS estimate without puffer transformation once
the set of candidate models with given sparsity pattern has been obtained. In
our test cases in chapter \ref{Results}, we only calculate BIC for models
without the puffer transformation and do hence not provide any formal tests for
this procedure.


\section{Numerical Solvers}
\label{chap_solvers}
To showcase the performance of SITE, we consider three example PDEs with
corresponding numerical solvers for which MDEs are to be identified from solver
output data, namely the advection, Burgers' and KdV equation. We performed all
simulations on a domain $\Omega: x \in [0,1]$ with periodic BCs $u(0,t)=u(1,t)$. 

\subsection{Advection Equation}
We consider the linear advection equation  \eqref{Advection_equation} for
$a = 1$ with a FTBS discretization scheme \eqref{FTBS_Scheme_theory}. The
linear accuracy requirement is well known to be
%
\begin{equation}
\mathrm{CFL} = \frac{a \Delta t}{\Delta x} \leq 1.
\end{equation}

\subsection{Burgers' Equation}
Burgers' equation is an important model equation in several fields, including
fluid dynamics and traffic flow calculations. Due to its similarity to the
Navier-Stokes equation and its tendency to develop shock solutions, it is
frequently considered as a test case for numerical algorithms in the
literature. To demonstrate applicability of SITE to nonlinear equations
and more advanced discretization schemes, we discretize the inviscid Burgers'
equation
%
\begin{equation}
    \label{Burgers_equation}
    u_t + \left(\dfrac{u^2}{2}\right)_x = 0
\end{equation}
by a MacCormack predictor-corrector scheme \cite{MACCORMACK.1969}, which is
second order accurate in space and time:
%
\begin{align}
    \begin{split}
        \tilde{u}_i^{j+1} = & u_i^j - h \left(\frac{(u^j_{i+1})^2}{2} - \frac{(u^j_{i})^2}{2} \right), \\
        u_i^{j+1} = & u_i^j - \frac{h}{2} \left[\left(\frac{(u^j_{i+1})^2}{2} - \frac{(u^j_{i})^2}{2} \right) + \left(\frac{(\tilde{u}^j_{i})^2}{2} - \frac{(\tilde{u}^j_{i-1})^2}{2} \right) \right].
        \label{Burgers Mac Cormack}
        \end{split}
\end{align}
The linear stability criterion is given by 
%
\begin{equation}
\mathrm{CFL} = \frac{|u|_{max} \Delta t}{\Delta x} \leq 1.
\end{equation}
For derivation of MDEs, the predictor-corrector scheme
\eqref{Burgers Mac Cormack} is rewritten into a single equation
\cite{Lerat.1974}
%
\begin{align}
    \begin{split}
        0 &  =  \frac{u_i^{j+1} - u_i^j}{\Delta t} + \frac{(u_{i+1}^j)^2 - (u_{i-1}^j)^2}{4 \Delta x} \\
            & - \frac{\Delta t}{2} \left( \frac{u_i^j + u_{i-1}^j}{2} \frac{(u_{i+1}^j)^2 - 2 (u_{i}^j)^2 + (u_{i-1}^j)^2}{2 \Delta x^2} + \frac{u_{i}^j - u_{i-1}^j}{\Delta x} \frac{(u_{i+1}^j)^2 - (u_{i-1}^j)^2}{4 \Delta x} \right) \\
            & + \frac{\Delta t^2}{2} \frac{(u_{i+1}^j)^2 - (u_{i-1}^j)^2}{4 \Delta x} \frac{(u_{i+1}^j)^2 - 2 (u_{i}^j)^2 + (u_{i-1}^j)^2}{2 \Delta x^2}.
    \end{split}
\end{align}

\subsection{Korteweg-de Vries Equation}
The KdV equation describes the asymptotic behavior of one-dimensional waves
with small amplitudes. The modeled physical phenomena include shallow water
waves and magneto-hydrodynamic waves in a plasma \cite{Zabusky.1965}. This test
case demonstrates the ability of SITE to identify first MDEs, whose higher
order temporal derivatives are not substituted with spatial derivatives. We
discretize the KdV equation
%
\begin{equation}
    u_t + 6 u u_x + u_{xxx} = 0
\end{equation}
using the Zabusky and Kruskal scheme \cite{Zabusky.1965}, which is second order
accurate in space and time:
%
\begin{align}
    \begin{split}
        \label{eq_KdV_scheme}
        u_i^{j+1} &= u_i^{j-1} - 2 h (u_{i+1}^j + u_i^j + u_{i-1}^j) (u_{i+1}^j - u_{i-1}^j) \\
        & - \frac{h}{\Delta x^2}  (u_{i+2}^j - 2u_{i+1}^j + 2 u_{i-1}^j - u_{i-2}^j).
    \end{split}
\end{align}
Due to the central in time approximation of \eqref{eq_KdV_scheme}, an
uncentered time discretization has to be used for the first time step \cite{Taha.1984}.
%
\begin{align}
    \begin{split}
        u_i^{1} &= u_i^{0} - h (u_{i+1}^0 + u_i^0 + u_{i-1}^0) (u_{i+1}^0 - u_{i-1}^0) \\
        & - \frac{h}{2 \Delta x^2}  (u_{i+2}^0 - 2u_{i+1}^0 + 2 u_{i-1}^0 - u_{i-2}^0)
    \end{split}
\end{align}
The linear stability criterion corresponding to the Zabusky and Kruskal scheme
is much more restrictive,
%
\begin{equation}
    \frac{\Delta t}{\Delta x} |-2 u_{\mathrm{max}} + \frac{1}{\Delta x^2} | \leq \frac{2}{3\sqrt{3}},
    \label{eq:KdV_stability_criterion}
\end{equation}
thus $\Delta t$ scaling with $\Delta x^3$ \cite{Taha.1984}.

\subsection{Verification}
We utilized the method of manufactured solutions \cite{Salari.2000} for
verification of the custom solvers to eliminate the possibility of incorrect
results due to a wrong solver implementation. For all solvers, we chose the
manufactured solution $\overline{u}(x,t) = \sin(2 \pi (x+t)) + 0.001$. We
refined $\Delta x$ and $\Delta t$ together with a constant $\mathrm{CFL}=0.1$
for the advection and Burgers' equation and $\mathrm{CFL}=10^{-10}$ for the
KdV equation, due to its more restrictive stability requirement. The solver
convergence plots are shown in fig. \ref{MMS_Convergence_Solvers} in the
Appendix confirming accuracy orders in space and time of 1, 2 and 2 for FTBS, 
MacCormack and the Zabusky and Kruskal scheme, respectively.


\section{Numerical Test Cases}
\label{Results}
We investigate the properties of SITE based on the three test cases outlined in
chapter \ref{chap_solvers}. The analysis focuses on the default setting of
SITE, which uses FoBa with spline initialization and without puffer
transformation.  All test cases exhibit the same structure: First, we outline
the discretization parameters of the data generating simulation as well as the
design of the library $\bm \Theta(\mathbf{u})$. Second, the MDE predicted by
the default setup of SITE and its accuracy is presented in tabular form. Next,
we study the impact of preconditioning steps and its interplay with the choice
of sparse regression algorithm. We conclude each test case with an analysis on
the impact of the simulation grid on regression accuracy and on BIC model
selection for the default setting. BIC model selection is compared to an
optimal procedure selecting always the optimal model from the set proposed by
FoBa. We defined the optimal model as the model with maximum number of correct
terms while not introducing any incorrect terms.

\subsection{Algorithmic Implementation Details}
\begin{figure}[t]
    \centering
    \scalebox{0.8}{\input{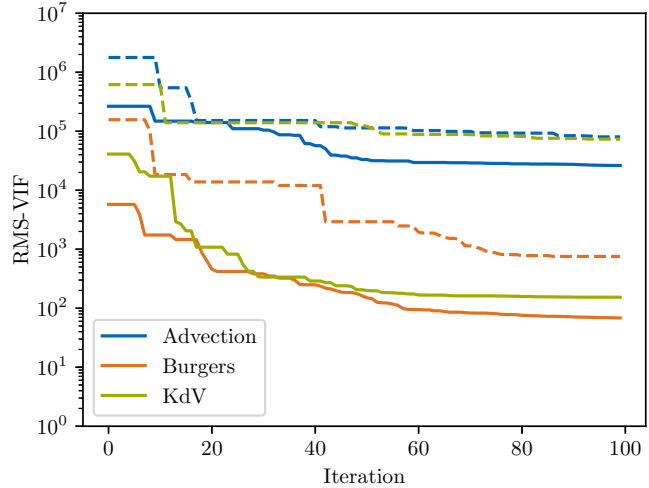}}
	\caption{NURBS particle swarm optimization for the three test cases. Solid
			 lines correspond to the training sets and dashed lines to the test
			 sets.}
    \label{VIF_spline_optimization}
\end{figure}

This section summarizes key implementation details of the following test cases.
We use \textit{Python 3.6} \cite{team2015python} with double precision numbers 
for all computations. In all our test cases we have $n > p$, but the method
extends in principle to $n < p$. Note however that we provide no formal tests
for this case. Derivatives are approximated by 8th order accurate finite
difference stencils from the \textit{findiff} \cite{Baer.2018} package. We only
include grid-points for which centered stencils are available. Non-centered
stencils introduce additional error in the derivative approximations, hence
impairing our ability to find the maximal number of candidate terms. The spline
initialization uses 8th order NURBS with 15 knots within $\Omega$ for the
trainingset and 11 knots for the testset. We optimized NURBS with the particle
swarm optimization algorithm implemented in the \textit{pyswarms}
\cite{JamesV.Miranda.2018} package with 50 particles, 100 iterations and
default optimization parameters. Fig. \ref{VIF_spline_optimization} illustrates
the VIF optimization progress for our test cases.

For comparison, a Gaussian bell curve is considered
%
\begin{equation}
	\label{Gauss_Init}
	u(x,0) = e^{-50(x-0.5)^2} + e^{-50(x+0.5)^2} + e^{-50(x-1.5)^2} \ ;  \quad x \in [0,1],
\end{equation}
which unlike $\sin(kx)$ from the example in section \ref{sec:multicolinearity}
is not inherently collinear to its derivatives.\\

We compare the maximum performance of the four algorithms from section
\ref{Theo Sparse Regression}. To avoid any bias resulting from model selection,
we iterate over their respective hyperparameters and evaluate the accuracy of
the predicted set of models. For SR3, we chose $L_0$ regularization. We stopped
the optimization process for Lasso and SR3 after 10000 iterations to keep
computational effort within reasonable bounds.

To assess the prediction accuracy of the sparse regression algorithms
considered here, we distinguish between correct terms, which are included in
the respective analytically derived MDE, and incorrect terms, which are not
included. If no incorrect term is included in the predicted model
$\mathcal{M}$, mean absolute error (MAE) and mean relative error (MRE) are
calculated from the analytically derived MDE weights $\bm{\xi}$,
%
\begin{equation}
	\mathrm{MAE} = \frac{1}{p_{\mathcal{M}}} \sum_{i=1}^{p_\mathcal{M}} |\xi_{i,\mathcal{M}} - \xi_i| \ ; \quad \mathrm{MRE} = \frac{1}{p_{\mathcal{M}}} \sum_{i=1}^{p_\mathcal{M}} \left|\frac{\xi_{i,\mathcal{M}} - \xi_i}{\xi_i} \right|.
\end{equation}
Note that leading order truncation error terms dominate MAE, while the highest
order terms dominate MRE (e.g. table \ref{tab:Advection_model}). Due to the
fact that deviations of even orders of magnitudes in the smallest term would be
invisible when adhering to MAE, we deem MRE to be more appropriate to judge the
accuracy of sparse regression algorithms. For practitioners who are primarily
interested in the leading order truncation error, MAE does however offer
valuable information.

We calculate the empirical order of identified terms from
$\xi_m \sim \mathcal{O}(\Delta x^k)$ using predictions from two simulations
with spatial width $\Delta x_1$ and $\Delta x_2$,
\begin{equation}
k = \frac{\log(\xi_{m1}/\xi_{m2})}{\log(\Delta x_1 / \Delta x_2)}.
\label{eq_p_results}
\end{equation}

\subsection{Advection Equation with FTBS}
To demonstrate the limits of accuracy of SITE, we study the linear advection
equation \eqref{Advection_equation} with a FTBS discretization scheme
\eqref{FTBS_Scheme_theory} with the objective to identify its third MDE
\eqref{Modified_Advection_equation}. The discretization parameters are
$\mathrm{CFL}=0.01$ on a grid $(n^x, n^t) = (300, 17)$, yielding 5 time steps
after data padding. We build two libraries: a small one contains $u$ and all
its spatial derivatives up to order 6. The large one appends all of those
combinations of derivatives that add up to a given cumulative order for all
cumulative orders up to 6, e.g. augmenting the library by $u_x^3$ and
$u_x u_{xx}$ for a cumulative order of 3. Next, these basis functions are
multiplied by $u^k$ for $k$ up to 6 and an intercept is added. This yields
$p=49$ candidate terms in the small library and $p=210$ in the large one.
Note that the maximum number of correct terms included in the libraries is 6.
The spline is optimized with respect to the small library. For the large
library, $n_{\mathrm{eff}} \approx p$, presuming our coarse estimate
$n_{\mathrm{eff}} = n^x$ and considering intra time step correlation.\\

Table \ref{tab:Advection_model} summarizes the predicted MDE of SITE in the
default setting for the large library case, demonstrating highly accurate
predictions. The empirical order is calculated pairwise from eq.
\eqref{eq_p_results} with a sequence of $n^x = (200, 300, 400, 500)$. The
obtained orders are averaged afterwards.

\begin{table}[t]
	\centering
	\caption{Summary of SITE default setup prediction for the advection equation}
	\label{tab:Advection_model}
	\begin{tabular}{lcccc}
		\toprule
		& analytical weight & absolute error & relative error & empirical order \\
		\midrule
		$v_x$		 & $-1$ & $7.89 \cdot 10^{-12}$ & $7.89 \cdot 10^{-12}$ & $0.00$ \\
		$v_{xx}$	 & $1.65 \cdot 10^{-3}$ & $8.55 \cdot 10^{-14}$ & $5.18 \cdot 10^{-11}$ & $1.00$ \\
		$v_{xxx}$	 & $-1.80 \cdot 10^{-6}$ & $6.48 \cdot 10^{-14}$ & $3.61 \cdot 10^{-8}$ & $2.00$ \\
		$v_{xxxx}$	 & $1.44 \cdot 10^{-9}$ & $2.50 \cdot 10^{-16}$ & $1.74 \cdot 10^{-7}$ & $3.00$ \\
		$v_{xxxxx}$	 & $-8.80 \cdot 10^{-13}$ & $1.57 \cdot 10^{-16}$ & $1.79 \cdot 10^{-7}$ & $4.00$ \\
		$v_{xxxxxx}$ & $4.04 \cdot 10^{-16}$ & $1.59 \cdot 10^{-19}$ & $3.94 \cdot 10^{-4}$ & $5.00$ \\
		\bottomrule
	\end{tabular}
\end{table}

\begin{figure}[t]
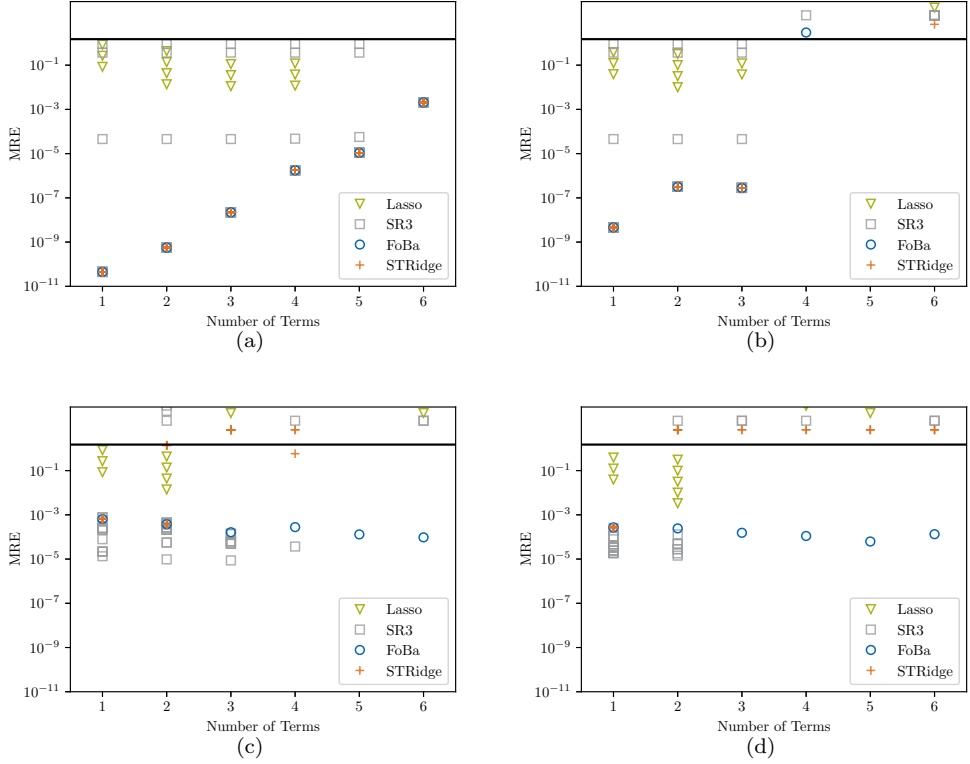

	\captionsetup[subfloat]{farskip=2pt,captionskip=1pt}
	\centering
	\subfloat[][]{ \scalebox{0.55}{
	\input{Regression_accuracy_comparison_Advection_norm_precondition_spline_small_Bib.pgf}}}
	\subfloat[][]{ \scalebox{0.55}{
	\input{Regression_accuracy_comparison_Advection_norm_precondition_gauss_small_Bib.pgf}}} \\
	\subfloat[][]{ \scalebox{0.55}{
	\input{Regression_accuracy_comparison_Advection_norm_spline_small_lib.pgf}}}
	\subfloat[][]{ \scalebox{0.55}{
	\input{Regression_accuracy_comparison_Advection_norm_Gauss_small_LIb.pgf}}}
	\caption{MRE of sparse regression algorithms as a function of the
			 number of terms included in a model from an iteration over the
			 respective hyperparameters for the advection equation and the
			 small library. Models above the black horizontal line
			 ($y = 10^{0}$) contain at least one incorrect term and are sorted
			 for visualization purposes. Their respective y-value has no
			 quantitative meaning. The considered setups include spline
			 initialization with puffer transformation (a), Gauss
			 initialization with puffer transformation (b), spline
			 initialization without puffer transformation (c) and Gauss
			 initialization without puffer transformation (d).}
	\label{fig:accuracy_advection}
\end{figure}

Fig. \ref{fig:accuracy_advection} compares the considered sparse regression
algorithms with respect to MRE and the number of correctly identified terms in
the small library case. Lasso consistently predicts less accurate models
compared to the other sparse regression algorithms due to the bias from the
$L_1$ regularization. Puffer transformation together with spline initialization
(a) improves MRE and the number of identified terms for all algorithms. The
only exception is FoBa with 6 terms in the model, exhibiting slightly reduced
MRE. This might be caused by noise inflation from the puffer transformation
that impacts terms with small signal-to-noise ratio the most. If STRidge
identifies the same sparsity pattern as FoBa, their results are identical, as
both rely on OLS to predict $\bm \xi$. Interestingly, when using puffer, FoBa,
STRidge and SR3 yield identical peak accuracy across the range of terms in
models. SR3 additionally proposes less accurate models due to the extra degree
of freedom from the hyperparameter $\gamma$. Without puffer (c), SR3 is more
accurate than FoBa when comparing models with a small number of terms, most
likely due to its built-in preconditioning. In contrast, FoBa proposes a model
that includes all 6 correct terms, whereas the maximum number of correct terms
predicted by SR3 is 4. spline initialization yields similar results to the
Gauss initial condition \eqref{Gauss_Init} if no puffer transformation is
applied (d). In conjunction with puffer transformations however (b), spline
initialization significantly improves regression results. This behavior stems
from smaller multicollinearity that induces less noise by the puffer
transformation.

\begin{figure}[H]
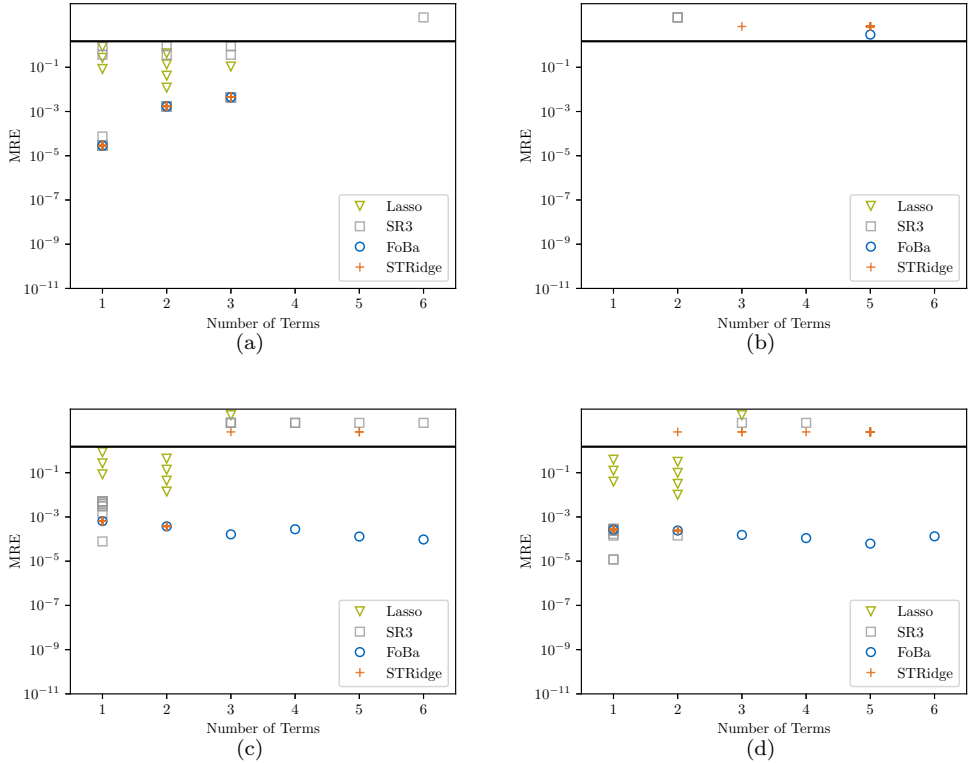

	\captionsetup[subfloat]{farskip=2pt,captionskip=1pt}
    \centering
    \subfloat[][]{ \scalebox{0.55}{
    \input{Regression_accuracy_comparison_Advection_norm_precondition_Spline_Large_Bib.pgf}}}
    \subfloat[][]{ \scalebox{0.55}{
    \input{Regression_accuracy_comparison_Advection_norm_precondition_Gauss_large_Bib.pgf}}}\\
    \subfloat[][]{ \scalebox{0.55}{
    \input{Regression_accuracy_comparison_Advection_norm_Spline_large_Lib.pgf}}}
    \subfloat[][]{ \scalebox{0.55}{
    \input{Regression_accuracy_comparison_Advection_norm_Gauss_large_Bib.pgf}}}
	\caption{MRE of sparse regression algorithms as a function of the number of
			 terms included in a model from an iteration over the respective
			 hyperparameters for the advection equation and the large library.
			 Models above the black horizontal line ($y = 10^{0}$) contain at
			 least one incorrect term and are sorted for visualization
			 purposes. Their respective y-value has no quantitative meaning.
			 The considered setups include spline initialization with puffer
			 transformation (a), Gauss initialization with puffer
			 transformation (b), spline initialization without puffer
			 transformation (c) and Gauss initialization without puffer
			 transformation (d).}
    \label{fig:accuracy_advection_large_Bib}
\end{figure}

When employing the large library (fig. \ref{fig:accuracy_advection_large_Bib}),
application of the puffer transformation significantly degrades regression
performance. This applies in particular if no spline initialization is used,
where none of the algorithms can identify a single correct term. We attribute
this behavior to the error inflating property of puffer, which becomes more
dominant for increasing VIF due to the additional candidate terms. Larger
multicollinearity from the Gauss initialization aggravates this effect. Without
a puffer transformation, the results are similar to the small library case,
with a slight decline in performance for SR3.

Fig. \ref{fig_Advection_resolution}(a) shows the impact of variations in $n^x$
on regression accuracy and on the ability of BIC in selecting the optimal model
for 5 time steps in the default setup. The proposed set of models always
contains the model with the maximum number of correct terms, however BIC fails
to select it for $n^x > 600$. MRE decreases until $n^x = 400$ due to the
decreasing importance of higher order terms not included in the model. These
decrease faster in magnitude than the sought terms as the largest correct term
not in the library is $\mathcal{O}(\Delta x^6)$. For $n^x > 400$, the benefits
of this effect are outweighed by increased noise from the round-off error of
the finite difference approximations, resulting in increasing MRE. MAE, which
is dominated by large magnitude terms with high signal-to-noise ratio, is more
robust to this additional noise. The results for 100 time steps are comparable
to the 5 time step case (figure \ref{fig_Advection_resolution}(b)) with the
only noticeable difference being that BIC identifies the correct model up to
$n^x = 900$.

\begin{figure}[H]
	\centering
	\subfloat[][]{ \scalebox{0.52}{
		\input{Resolution_properties_Advection_norm17.pgf}
	}}
	\subfloat[][]{ \scalebox{0.52}{
		\input{Resolution_properties_Advection_norm112.pgf}
	}}
	\caption{Number of terms selected by BIC and by optimal choice from the
			 model candidates proposed by FoBa as a function of resolution for
			 the advection equation and 5 time steps (a) or 100 time steps (b).
			 Markers of models only containing correct terms are plotted in
			 blue, models containing at least one incorrect term are plotted in
			 red. BIC is represented by a cross ($\times$) and the optimal
			 choice by a circle ($\circ$). MRE and MAE are displayed for the
			 optimal model on the right axis.}
	\label{fig_Advection_resolution}
\end{figure}

\subsection{Burgers' Equation with MacCormack}
After inserting Taylor series into the MacCormack scheme
\eqref{Burgers Mac Cormack} for the Burgers' equation \eqref{Burgers_equation}
and substituting higher order temporal derivatives with spatial derivatives,
the third MDE is obtained \cite{Lerat.1974}.
%
\begin{equation}
	\label{MDE_Burgers}
	v_t + \left(\frac{v^2}{2} \right)_x - \Delta x^2 \left(v_{xxx}\frac{v}{6} (h^2v^2 - 1) + \frac{v_x v_{xx}}{2} (2 h^2 v^2 - hv - 1) + v_x^3\frac{h}{4} (2hv - 1) \right) + \mathcal{O}(\Delta x^3) = 0
\end{equation}
The discretization parameters are $\mathrm{CFL}=0.5$ on a grid
$(n^x, n^t) = (10000, 17)$, yielding 5 time steps after data padding. The
larger $n^x$ and CFL number result in a time step size in the same order of
magnitude as in the advection case. The library contains $u$ and all its
spatial derivatives up to order 3 as well as all of those combinations of
derivatives which add up to a given cumulative order for all cumulative orders
up to 3. These basis functions are then multiplied by $u^k$ for $k$ up to 3 and
an intercept is added, yielding $p=28$ candidate terms. Since
$\bm \Theta(\mathbf{u})$ can represent all second order truncation error terms
in the MDE \eqref{MDE_Burgers}, the maximum number of correctly identifiable
terms is 8.\\

Table \ref{tab:Burgers_model} summarizes the predicted MDE of SITE in the
default setting, demonstrating consistently accurate predictions of the
truncation error terms with relative errors in the order $10^{-4}$. The
empirical order is calculated by a sequence of
$n^x = (6000, 8000, 10000, 12000)$.

\begin{table}[t]
	\centering
	\caption{Summary of SITE default setup prediction for the Burgers' equation}
	\label{tab:Burgers_model}
	\begin{tabular}{lcccc}
		\toprule
		& analytical weight & absolute error & relative error & empirical order \\
		\midrule
		$v v_x$ 			& $-1$ & $5.71\cdot 10^{-10}$ & $5.71\cdot 10^{-10}$ & $0.00$ \\
		$v^3 v_{xxx}$ 		& $9.87\cdot 10^{-9}$ & $-9.36\cdot 10^{-12}$ & $9.49\cdot 10^{-4}$ & $2.01$ \\
		$v v_{xxx}$			& $-1.67\cdot 10^{-9}$ & $2.21\cdot 10^{-13}$ & $1.33\cdot 10^{-4}$ & $2.00$ \\
		$v^2 v_x v_{xx}$	& $5.92\cdot 10^{-8}$ & $3.96\cdot 10^{-11}$ & $6.69\cdot 10^{-4}$ & $2.00$ \\
		$v v_x v_{xx}$		& $-1.22\cdot 10^{-8}$ & $4.79\cdot 10^{-12}$ & $3.94\cdot 10^{-4}$ & $2.00$ \\
		$v_x v_{xx}$		& $-5.00\cdot 10^{-9}$ & $-7.40\cdot 10^{-13}$ & $1.48\cdot 10^{-4}$ & $2.00$ \\
		$v v_{x}^3$			& $2.96\cdot 10^{-8}$ & $-7.76\cdot 10^{-11}$ & $2.62\cdot 10^{-3}$  & $1.99$ \\
		$ v^3_{x}$			& $-6.08\cdot 10^{-9}$ & $1.74\cdot 10^{-12}$ & $2.86\cdot 10^{-4}$ & $2.00$ \\
		\bottomrule
	\end{tabular}
\end{table}

\begin{figure}[H]
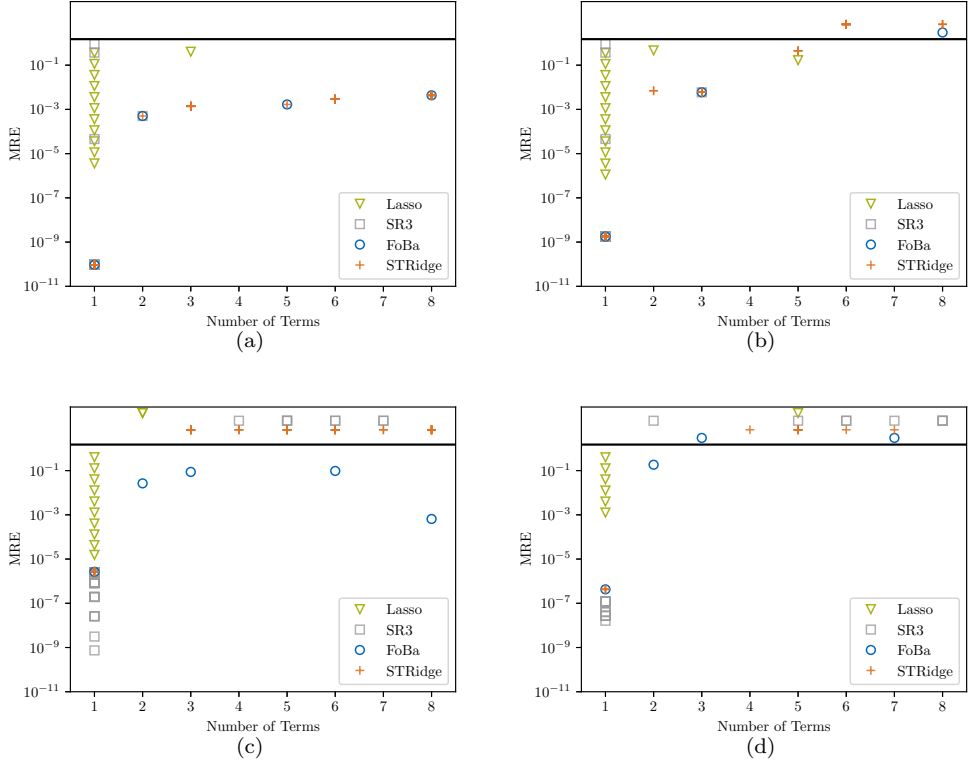

	\captionsetup[subfloat]{farskip=2pt,captionskip=1pt}
	\centering
	\subfloat[][]{ \scalebox{0.55}{
	\input{Regression_accuracy_comparison_Burgers_norm_precondition_spline.pgf}}}
	\subfloat[][]{ \scalebox{0.55}{
	\input{Regression_accuracy_comparison_Burgers_norm_precondition_Gauss.pgf}}} \\
	\subfloat[][]{ \scalebox{0.55}{
	\input{Regression_accuracy_comparison_Burgers_norm_spline.pgf}}}
	\subfloat[][]{ \scalebox{0.55}{
	\input{Regression_accuracy_comparison_Burgers_norm_Gauss.pgf}}}
	\caption{MRE of sparse regression algorithms as a function of the number of
			 terms included in a model from an iteration over the respective
			 hyperparameters for the Burgers' equation. Models above the black
			 horizontal line ($y = 10^{0}$) contain at least one incorrect term
			 and are sorted for visualization purposes. Their respective
			 y-value has no quantitative meaning. The considered setups include
			 spline initialization with puffer transformation (a), Gauss
			 initialization with puffer transformation (b), spline
			 initialization without puffer transformation (c) and Gauss
			 initialization without puffer transformation (d).}
	\label{fig:accuracy_burgers}
\end{figure}

Fig. \ref{fig:accuracy_burgers} compares the accuracy of sparse regression for
Burgers' equation. STRidge profits substantially from puffer transformation,
then being able to identify all terms correctly using Spline initialization and
outperforming FoBa with Gauss initialization. Without puffer, FoBa still
detects all terms using spline initialization, while all other algorithms can
only detect the term from Burgers' equation, but no truncation error terms.
Spline initialization considerably improves results with and without puffer
transformation. Interpretations are analogous to the advection test case.

\begin{figure}[H]
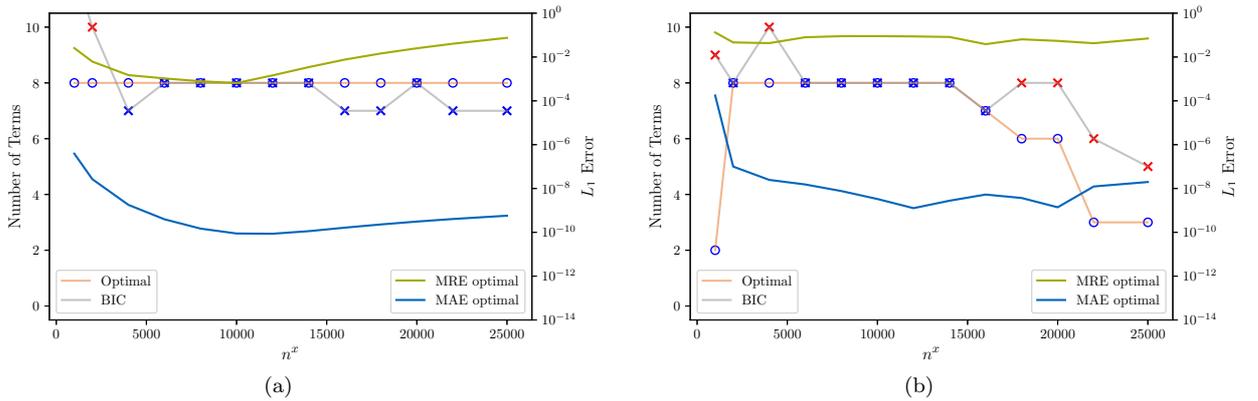

	\centering
	\subfloat[][]{ \scalebox{0.52}{
		\input{Resolution_properties_Burgers_norm17.pgf}
	}}
	\subfloat[][]{ \scalebox{0.52}{
		\input{Resolution_properties_Burgers_norm112.pgf}
	}}
	\caption{Number of terms selected by BIC and by optimal choice from the
			 model candidates proposed by FoBa as a function of resolution for
			 the Burgers' equation and 5 time steps (a) or 100 time steps (b).
			 Markers of models only containing correct terms are plotted in
			 blue, models containing at least one incorrect term are plotted in
			 red. BIC is represented by a cross ($\times$) and the optimal
			 choice by a circle ($\circ$). MRE and MAE are displayed for the
			 optimal model on the right axis.}
			 \label{fig_Burgers_resolution}
\end{figure}

Fig. \ref{fig_Burgers_resolution}(a) shows the impact of variations in $n^x$
with respect to regression accuracy and the model selection capabilities of BIC
for 5 time steps in the default setup. Similarly to the advection case, the
model with the maximum number of correct terms is always included in the
proposed set. However, BIC can only identify it within a range of
$n^x \in [6000, 14000] $, where the optimal error is the smallest. Given that
all truncation error terms are of the same order $\mathcal{O}(h^2)$, the
increased noise impacts not only MRE, but MAE as well. Since the highest
derivative to be approximated is 3, $n^x$ can be chosen much larger than in the
advection case before round-off error becomes dominant. Therefore, correct
higher order truncation error terms not included in the library can be driven
towards 0 effectively.

The results based on the simulation with 100 time steps is shown in fig.
\ref{fig_Burgers_resolution}(b). In contrast to previous examples, the model
with the maximum number of correct terms is often not included in the set of
models provided by FoBa outside the range $n^x \in [6000, 14000] $. MRE and MAE
increase significantly compared to the 5 time step case, indicating increased
noise in the data obtained from later stages of the simulation. Note that both
curves lose some meaning in areas where the number of terms in the optimal
model changes.

\subsection{Korteweg-de Vries Equation with Zabusky and Kruskal}

Inserting Taylor series expansions into the Zabusky and Kruskal
\cite{Zabusky.1965} discretization scheme \eqref{eq_KdV_scheme} yields the
first MDE.
%
\begin{align}
	\label{MDE_Kdv}
	\begin{split}
	& v_t + 6v v_x + v_{xxx} + \Delta x^2 \left(\frac{h^2}{6} v_{ttt} + \frac{1}{4}u_{xxxxx} + v v_{xxx} + 2 v_x v_{xx}\right) \\
	&+ \Delta x^4 \left(\frac{h^4}{120} v_{ttttt} + \frac{1}{40} v_{xxxxxxx} + \frac{1}{3} v_{xx}v_{xxx} + \frac{1}{6} v_x v_{xxxx} + \frac{1}{20} v v_{xxxxx}\right) + \mathcal{O}(\Delta x^6) = 0
	\end{split}
\end{align}
The discretization parameters are $\mathrm{CFL}=10^{-6}$ on a grid
$(n^x, n^t) = (100, 19)$, yielding 5 time steps after data padding. The
library contains $u$ and all its spatial derivatives up to order 7 as well as
all of those combinations of derivatives that add up to a given cumulative
order for all cumulative orders up to 3. These basis functions are then
multiplied by $u^k$ for $k$ up to 3 and an intercept is added. Our goal is to
represent a first MDE, thus 2nd and 3rd order time derivatives are appended,
yielding $p=68$ candidate terms. Note the increased padding width due to the
higher order time derivatives. $\bm \Theta(\mathbf{u})$ can represent 8 terms
from the first MDE up to 4th order \eqref{MDE_Kdv}
and one additional term from the 6th order truncation error, not shown in eq.
\eqref{MDE_Kdv}. We constructed $\bm \Theta(\mathbf{u})$ this way to
prove that first MDEs can be found with SITE. When including all terms of
$\mathcal{O}(h^4)$, we found all spatial derivatives, but could not identify
$v_{ttt}$, which is very small (table \ref{tab:KdV_model}) due to the small CFL
number enforced by the solver stability criterion.\\

\begin{table}[t]
	\centering
	\caption{Summary of SITE default setup prediction for the KdV equation}
	\label{tab:KdV_model}
	\begin{tabular}{lcccc}
		\toprule
		& analytical weight & absolute error & relative error & empirical order \\
		\midrule
		$v v_x$ 		& $-6$ & $2.25 \cdot 10^{-2}$ & $3.75 \cdot 10^{-3}$ & $-0.01$ \\
		$v_{xxx}$ 		& $-1$ & $-5.08 \cdot 10^{-6}$ & $5.08 \cdot 10^{-6}$ & $0.00$ \\
		$v_{ttt}$ 		& $-2.09 \cdot 10^{-15}$ & $7.24 \cdot 10^{-16}$ & $3.47 \cdot 10^{-1}$ & $2.03$ \\
		$v_{xxxxx}$ 	& $-2.5 \cdot 10^{-5}$ & $-7.47 \cdot 10^{-9}$ & $2.99 \cdot 10^{-4}$ & $2.00$ \\
		$v v_{xxx}$ 	& $-1 \cdot 10^{-4}$ & $1.85 \cdot 10^{-5}$ & $1.85 \cdot 10^{-1}$ & $1.96$ \\
		$v_x v_{xx}$ 	& $-2 \cdot 10^{-4}$ & $5.96 \cdot 10^{-6}$ & $2.98 \cdot 10^{-2}$ & $1.92$ \\
		$v_{xxxxxxx}$ 	& $-2.5 \cdot 10^{-10}$ & $-2.06 \cdot 10^{-12}$ & $8.26 \cdot 10^{-3}$ & $3.98$\\
		\bottomrule
	\end{tabular}
\end{table}

Table \ref{tab:KdV_model} summarizes the predicted MDE of SITE in the default
setting. The empirical order is calculated by a sequence of
$n^x = (88, 100, 112, 125)$. Note that regression accuracy is significantly
lower compared to both previous examples.

\begin{figure}[H]
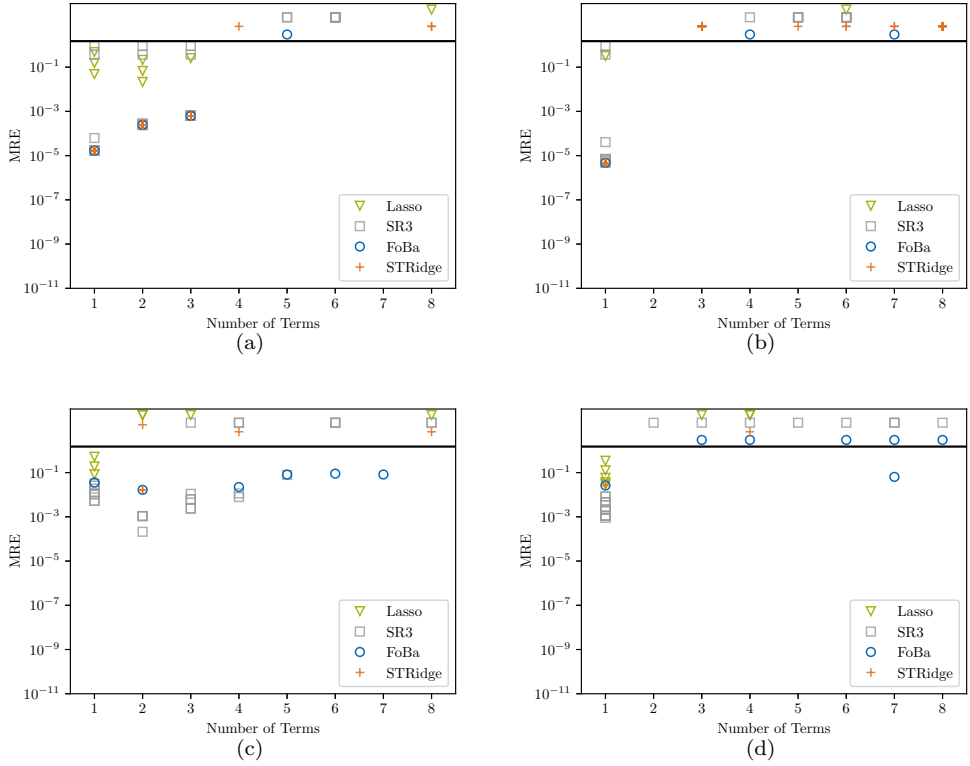

	\captionsetup[subfloat]{farskip=2pt,captionskip=1pt}
	\centering
	\subfloat[][]{ \scalebox{0.55}{
	\input{Regression_accuracy_comparison_KdV_norm_precondition_spline.pgf}}}
	\subfloat[][]{ \scalebox{0.55}{
	\input{Regression_accuracy_comparison_KdV_norm_precondition_gauss.pgf}}}\\
	\subfloat[][]{ \scalebox{0.55}{
	\input{Regression_accuracy_comparison_KdV_norm_spline.pgf}}}
	\subfloat[][]{ \scalebox{0.55}{
	\input{Regression_accuracy_comparison_KdV_norm_gauss.pgf}}}
	\caption{MRE of sparse regression algorithms as a function of the number of
			 terms included in a model from an iteration over the respective
			 hyperparameters for the KdV equation. Models above the black
			 horizontal line ($y = 10^{0}$) contain at least one incorrect term
			 and are sorted for visualization purposes. Their respective
			 y-value has no quantitative meaning. The considered setups include
			 spline initialization with puffer transformation (a), Gauss
			 initialization with puffer transformation (b), spline
			 initialization without puffer transformation (c) and Gauss
			 initialization without puffer transformation (d).}
	\label{fig:accuracy_kdv}
\end{figure}

The results of the sparse regression comparison in fig. \ref{fig:accuracy_kdv}
are mostly analogous to the other test cases. With applied puffer transformation, no 
algorithm can detect more than 3 terms.
This is significantly smaller than the number of terms identifiable without puffer. 
Gauss initialization yields poor results in this test case: the majority of predicted models 
only contains one of two terms of the KdV equation.

\begin{figure}[H]
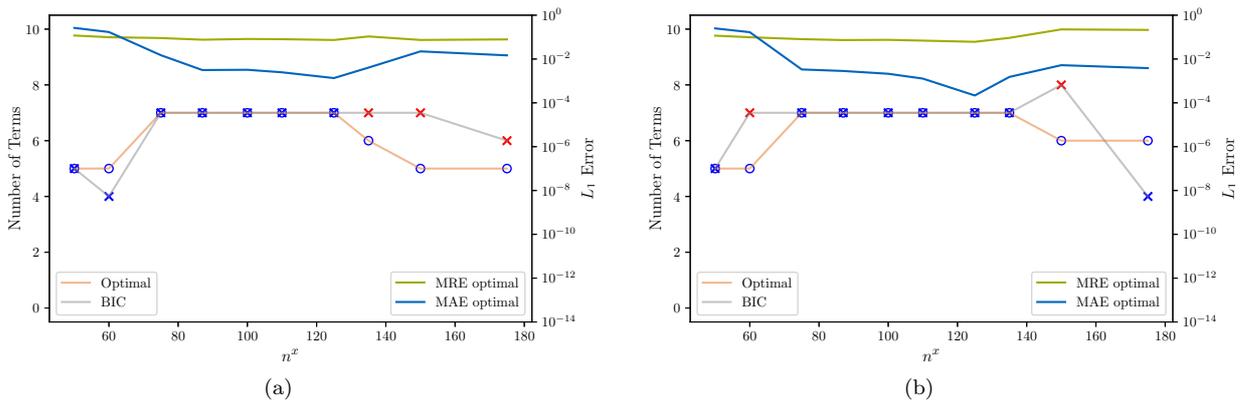

	\centering
	\subfloat[][]{ \scalebox{0.52}{
		\input{Resolution_properties_KdV_norm19.pgf}
	}}
	\subfloat[][]{ \scalebox{0.52}{
		\input{Resolution_properties_KdV_norm119.pgf}
	}}
	\caption{Number of terms selected by BIC and by optimal choice from the
			 model candidates proposed by FoBa as a function of resolution
			 for the KdV equation and 5 time steps (a) or 100 time steps (b).
			 Markers of models only containing correct terms are plotted in
			 blue, models containing at least one incorrect term are plotted in
			 red. BIC is represented by a cross ($\times$) and the optimal
			 choice by a circle ($\circ$). MRE and MAE are displayed for the
			 optimal model on the right axis.}
			 \label{fig_KdV_resolution}
\end{figure}

Fig. \ref{fig_KdV_resolution}(a) shows the dependency of SITE predictions on
spatial resolution in default setting for 5 time steps. The optimal resolution
range $n^x \in [75, 125]$ is very narrow and outside this range, the maximum
number of terms in models proposed by FoBa decreases instantly. Even for
$n^x \in [75, 125]$, one term within $\mathcal{O}(h^4)$ - which could be
represented by $\bm \Theta(\mathbf{u})$ - is not identified. Given that a 7th
order spatial derivative is to be approximated, $n^x $ is considerably limited
due to round-off error effects. Unlike in the Burgers case, driving higher
order terms towards 0 is therefore not possible, compromising the ability to
identify all terms within one order of magnitude. If we were not to include the
7th spatial derivative in the library, $n^x$ would still be limited due to the
restrictive stability criterion of the Zabusky and Kruskal scheme
\eqref{eq:KdV_stability_criterion}. As $\Delta t$ scales with $\Delta x^3$,
increasing $n^x$ would quickly decrease $\Delta t$ and therefore result in
increasing round-off errors for the time derivative approximations, for which
at least 3rd-order derivatives are to be approximated. Utilizing simulation
data from 100 time steps (fig. \ref{fig_KdV_resolution}(b)) instead of 5
slightly improves quality measures in this test case, including an extended
optimal resolution range up to $n^x = 137$.

\section{Discussion}
\label{chap_Discussion}
We discuss a few guidelines obtained from the advection, Burgers and KdV test
cases. The predictions of our SITE approach are determined by three factors:
Simulation data, choices in the algorithmic procedure and in the model
selection step. The simulation data has considerable influence, defining the
signal-to-noise ratio of the regression problem. $\Delta x$ and $\Delta t$
should be chosen to balance noise from higher order truncation error terms not
included in the library and noise due to round-off error of the finite
difference approximations. While we mainly focused on the  impact of $n^x$ in
chapter \ref{Results}, we found that prediction results strongly depend on a
sensible choice of $\Delta t$ as well. The range of $\Delta t$ and $\Delta x$
for maximum accuracy predictions differs significantly between problems. This
range is particularly narrow when high order spatial and temporal derivatives
are being approximated as in the KdV case. Round-off error substantially limits
the range of simulation parameters, curtailing the opportunity to drive high
order terms towards 0. Using higher precision numbers than double precision, as
it is possible with e.g. FORTRAN \cite{metcalf2004fortran} or Julia
\cite{bezanson2012julia}, has the potential for significant improvements in
both regression accuracy and number of identifiable terms. 

The number of simulation time steps impacts regression accuracy, the maximum
number of  correct terms proposed by FoBa and BIC model selection. However, the
extent of impact on these three quality measures differed considerably between
test cases and increasing the number of time steps did not reveal a clear
positive or negative trend. The mechanisms involved are not yet fully
understood and warrant further investigation. \\

The algorithmic setup is defined by appropriate choices of library,
preconditioning and sparse regression algorithm. In library construction, using
existing prior knowledge about the form of candidate functions is advisable.
Reducing the number of candidates usually leads to smaller multicollinearity,
which is beneficial to sparse regression. Performing an iterative approach by
constructing the library over multiple successive iterations might be helpful.
At this point, a word of caution is in order: In case the algorithm has
uncharacteristically poor performance this may be due to an erroneous
construction of the candidate term library. If a correct, strong impact term is
not included in the library, we have seen that this term is then often
approximated by a linear combination of incorrect terms from the library.

The test cases demonstrate that the benefits of preconditioning depend on the
problem considered. Puffer transformation can improve both regression accuracy
and the number of identified terms if the signal-to-noise ratio of terms
included in the model is high. However, the noise inflation property quickly
outweighs these benefits if multicollinearity is severe. Using spline
initialization together with puffer is essential to decrease the amount of
additional noise from multicollinearity. Generalizations of $\mathbf{F}$,
following Jia et al. \cite{Jia.2015}, might aid in the reduction of the noise
inflation effect. Spline initialization without puffer behaves the opposite
way. While there is only a minor benefit in low noise problems like the
advection case, the benefits are considerable in higher noise problems
(Burgers, KdV). If an unconstrained choice of $u(x,0)$ is infeasible due to
problem constraints, spline initialization can be skipped. 

Preconditioning choices are also linked to the selected sparse regression
algorithm. For practical applications, a robust default setup seems to be FoBa
with spline initialization and without the puffer transformation. Apart from
robust term detection, FoBa yields highly accurate predictions even without
puffer transformation. This is only surpassed by SR3 for models with a small
number of terms. FoBa facilitates model selection in comparison to Lasso or
SR3. While the former only proposes a few models, the latter yield different
models for every set of hyperparameters. If a different algorithm than FoBa is
used, puffer transformation often improves the results both in accuracy and
number of identified terms. For our test cases, we did not find any evidence
that STRidge was more reliable than FoBa
\footnote{We assume that the reduced reliablilty of FoBa with respect to
          STRidge reported by Rudy et al. \cite{Rudy.2017} was due to an
          incorrect implementation of FoBa. For our implementation, we
          therefore made adjustments according to the work of Zhang \cite{Zhang.2009}.}. \\

Model selection using information criteria has previously been introduced to
the PDE-FIND framework by Mangan et al. \cite{Mangan.2017}. Their approach is
based on a testset of $N \sim 100$ simulations $\{u^i(x,t)\}_{i=1}^N$ with
varying initial conditions $u^i(x,0)$. For each model from $\{M_j\}_{j=1}^K$,
the dynamics of model $j$ are integrated in time for each $u^i(x,0)$, yielding
$u^i_k(x,t)$. For each simulation in the testset, a residual between $u^i(x,t)$
and $u^i_k(x,t)$ is calculated and the RMS error of these residuals is used as
likelihood in the Akaike information criterion (AIC) \cite{Akaike.1974}. Since
each simulation is considered as one sample in AIC, a sufficiently large $N$
has to be chosen, resulting in large computational effort for time integration
of the candidate models. In order to avoid the computational effort and the
issues related to time integration of MDEs discussed in section
\ref{chap_MDEA}, we calculate BIC directly from the linear model
\eqref{SINDy System of equations} with respect to a single test simulation. A
major advantage of the procedure of Mangan et al. \cite{Mangan.2017} is that
rather uncorrelated samples can be used in AIC, hence eliminating the need to
estimate $n_{\mathrm{eff}}$, which becomes a free parameter of the model
selection step of SITE.

Estimating $n_{\mathrm{eff}}$ with $n^x$ is very coarse: $n_{\mathrm{eff}}$
clearly does not increase linearly with $n^x$ due to increased intra timestep
correlation. However, we found our results to be robust with respect to
$n_{\mathrm{eff}}$ as long as the order of magnitude was roughly correct. For
problems where convection is not the dominant effect, a different choice of
$n_{\mathrm{eff}}$ should be made on a problem specific basis. For alternatives
to the proposed rule of thumb, one might consider estimation methods from the
literature e.g. \cite{Berger.2014}, or \cite{WEAKLIEM.1999}. There is
legitimate criticism aimed at the BIC due to its implicit dependence on a prior
distribution which can substantially deviate from the prior distribution a
considerate investigator would choose \cite{WEAKLIEM.1999}. However, we found
BIC to be capable of identifying the optimal model for a range of $n_x$ in all
test cases. Instead of BIC, AIC could potentially be a natural second choice,
yielding similar model selection from our experience. There exists a vast
amount of literature on improvements to the BIC, such as EBIC \cite{Chen.2008}
and MBIC \cite{WEAKLIEM.1999}, allowing for custom adjustments of the BIC to
the individual problem at hand.

\section{Conclusion}
\label{chap_Conclusion}
We presented SITE, a novel data-driven approach to modified differential
equation analysis. Its effectiveness in discovering first and third MDEs was
demonstrated in various test cases. The current implementation of SITE is by no
means optimal, as neither the preconditioning steps nor the model selection
procedure have been optimized. We still showed that high quality results can be
obtained, underlining the promise of the approach for extending the MDEA
toolbox to discretization schemes, where analytic derivation of MDEs is
infeasible. Applications that might benefit from SITE include optimization of
numerical discretization schemes and ILES turbulence modeling. A stepping stone
for the approach will be a deeper understanding of its limitations, glimpses of
which we were able to witness with Burgers' equation and the KdV equation.
Other research directions deemed worthy of inquiry are the handling of more
nonlinear problems, higher order derivatives, truncation errors in multiple
dimensions, application to the finite element method and application to the
finite volume method.


\section{Acknowledgements}

This project has received funding from the German Research Council (DFG)
under grant agreement No. 326472365. All data and codes used in this manuscript are publicly
available on Github at \newline \url{https://www.github.com/tumaer/truncationerror}.

\bibliography{Bibliography}
\bibliographystyle{siam}


\appendix

\section{}

\begin{figure}[H]
    \centering
    \scalebox{0.9}{\input{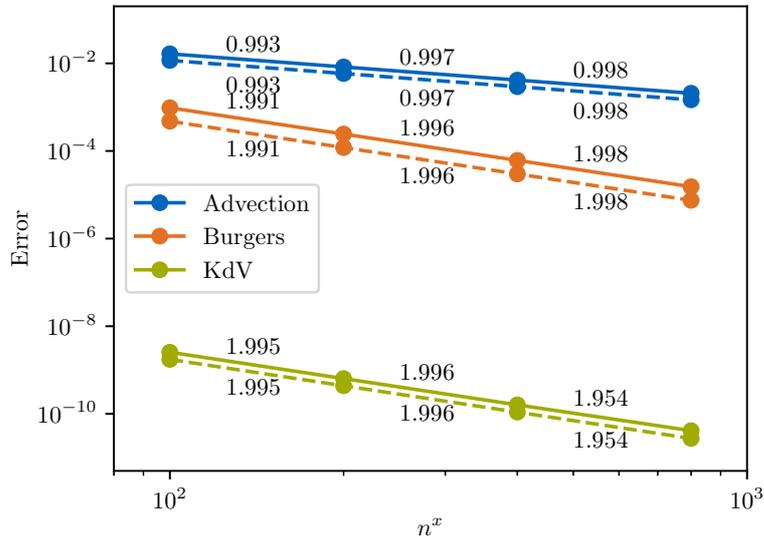}}
    \caption{$L_0$ (solid line) and $L_2$ error (dashed line)  as a function of grid resolution for MMS with empirical orders noted for each refinement step. $\mathrm{CFL}=0.1$ for FTBS and MacCormack and $\mathrm{CFL}=10^{-10}$ for the Zabusky and Kruskal scheme due to the restrictive stability criterion. The errors are computed at $t_{\mathrm{test}} = 0.1$ for FTMS and MacCormack and $t_{\mathrm{test}} = 10^{-8}$ for the Zabusky and Kruskal scheme}
    \label{MMS_Convergence_Solvers}
\end{figure}

\end{document}